\documentclass[12pt]{amsart}   	% use "amsart" instead of "article" for AMSLaTeX format
\usepackage{geometry}                		% See geometry.pdf to learn the layout options. There are lots.
\input{epsf.tex}
\usepackage{amsfonts,amssymb,verbatim,amsmath,amsthm,latexsym,textcomp,amscd}
\usepackage{latexsym,amsfonts,amssymb,epsfig,verbatim}
\usepackage{amsmath,amsthm,amssymb,latexsym,graphics,textcomp, graphicx,color}

\newtheorem{theorem}{Theorem}[section]
\newtheorem{thm}{Theorem}[section]

\newtheorem{defn}{Definition}[section]

\newtheorem{rem}[theorem]{Remark}

\newtheorem{notation}[theorem]{Notation}

\newtheorem{lemma}[theorem]{Lemma}

\newtheorem{proposition}[theorem]{Proposition}
\newtheorem{prop}[theorem]{Proposition}

\newtheorem{cor}[theorem]{Corollary}

\newtheorem{conjecture}[theorem]{Conjecture}

\newtheorem{example}[theorem]{Example}
\newtheorem{exe}[theorem]{Exercise}

\newtheorem{problem}[theorem]{Problem}
\newtheorem{question}[theorem]{Question}

\def\B{\mathbf B}
\def\k{\mathbf k}
\def\BB{\mathcal B}

\def\C{\mathbb C}
\def\P{\mathbb P}
\def\H{\mathbb H}

\def\Q{\mathbb Q}

\def\R{\mathbb R}
\def\Z{\mathbb Z}

\def\D{\partial}
\def\ol{\overline}

\def\<{\langle}
\def\>{\rangle}
\def\geo{\partial_\infty}

\def\F{{\mathcal F}}

\def\O{{\mathcal O}}

\def\al{\alpha} 
\def\Ga{\Gamma}
\def\Om{\Omega}
\def\ga{\gamma}
\def\La{\Lambda}
\def\la{\lambda}
\def\eps{\epsilon}
\def\Si{\Sigma}
\def\si{\sigma}

\def\embed{\hookrightarrow} 
\def\id{\operatorname{id}}
\def\inte{\operatorname{int}}
\def\hull{\operatorname{hull}}
\def\dim{\operatorname{dim}}
\def\card{\operatorname{card}}
\def\Isom{\operatorname{Isom}}
\def\arccosh{\hbox{arccosh}}
\def\cd{\operatorname{cd}}

\title{Lectures on complex hyperbolic Kleinian groups}
\author{Michael Kapovich}
\address{M.K.: Department of Mathematics, UC Davis, One Shields Avenue, Davis CA 95616, USA}
\email{kapovich@math.ucdavis.edu}
\thanks{During the work on this paper the  author was partly supported by the NSF grant  DMS-16-04241.} 

\date{\today}							% Activate to display a given date or no date

\begin{document}

%\section{}
%\subsection{}

\begin{abstract}
These are lectures on discrete groups of isometries of complex hyperbolic spaces, aimed to discuss interactions between the function theory on complex hyperbolic manifolds and the theory of discrete groups. 
\end{abstract}

\maketitle

\section{Introduction}

These notes are based on a series of lectures I gave at the workshop ``Progress in Several Complex Variables,'' held in KIAS, Seoul, Korea, in October of 2019. 
It is useful to read the notes in conjunction with my (longer) survey of discrete isometry groups of real hyperbolic spaces, \cite{Kapovich2008}, since most issues 
in the real and complex hyperbolic setting are quite similar. The theory of complex hyperbolic manifolds and complex hyperbolic Kleinian groups (aka discrete holomorphic isometry groups of complex hyperbolic spaces $\H^n_\C$) is a rich mixture of Riemannian and complex geometry, topology, dynamics, symplectic geometry and complex analysis. The choice of topics covered in these lectures is governed by my personal taste  and is, by no means, canonical: It is geared towards a discussion of interactions between the function theory on complex hyperbolic manifolds and the geometry/dynamics of 
complex hyperbolic Kleinian groups (sections \ref{sec:interactions} and \ref{sec:conjectures}). I refer the reader to \cite{CNS, Goldman, Gol, Epstein, McReynolds, Parker, Parker2, Paupert, Sch} for further discussion of geometry of complex hyperbolic spaces and their discrete isometry groups. The bibliography of complex hyperbolic Kleinian groups appearing at the end of these notes is long but is not meant to be exhaustive, my  apologies to everybody whose papers are omitted. 

\medskip 
%{\bf Acknowledgements.} During the work on this paper the  author was partly supported by the NSF grant  DMS-16-04241. 

\tableofcontents  

\section{Complex hyperbolic space}

Most of the basic material on geometry of complex hyperbolic spaces can be found in Goldman's book \cite{Gol}; I also refer the reader to \cite{Epstein, Parker, Paupert} for shorter introductions. 

\medskip 
Consider the vector space $V=\C^{n+1}$ equipped with the pseudo-hermitian bilinear form
$$
\<z, w\>=-z_0 \bar{w}_0 + \sum_{k=1}^n z_k \bar{w}_k.   
$$
Set $q(z):= \<z, z\>$. This quadratic form has signature $(n,1)$.  
Define the {\em negative light cone} $V_-:= \{z: q(z)<0\}$. Consider the complex projective space $\P^n:= P V$, the projectivization of $V$, and the projection $p: z\mapsto [z]\in \P^n$. The projection $\B^n:=p(V_-)$ is an open  ball in $\P^n$. In order to see this, consider the affine hyperplane in $\C^{n+1}$ given by $A=\{z_0=1\}$ (and equipped with the standard Euclidean hermitian metric). 
Then $V_-\cap A$ is the open unit ball in $A$ centered at the origin. This intersection projects diffeomorphically to $p(V_-)$.

The tangent space $T_{[z]} \P^n$ is naturally identified with $z^\perp$, the orthogonal complement of $\C z$ in $V$, taken with respect to $\<\cdot, \cdot\>$. 
If $z\in V_-$, then the restriction of $q$ to $z^\perp$ is positive-definite, hence, $\<\cdot, \cdot\>$ 
project to a hermitian metric $h$ (also denoted $\<\cdot, \cdot\>_h$) on $\B^n$. 
From now on, I will always equip $\B^n$ with the hermitian metric $h$ and let $d$ denote the corresponding distance function on $\B^n$. 

\begin{defn}
The {\em complex hyperbolic $n$-space} $\H^n_{\C}$ is $(\B^n,h)$. 
\end{defn}

I next describe the hermitian metric $h$ on $\B^n$ using the coordinates $(z_1,...,z_n)$ on $A$. First, regarding $\B^n$ as a subset of the affine hyperplane $A$, for a vector $y\in T_x \B^n$ 
we have
$$
\<y, y\>_h= \frac{\<x, x\> \<y, y\> - \<x, y\> \<y, x\>}{-\<x, x\>^2}. 
$$
Setting $x= (1, z)$, $z\in \C^n$, and denoting $u\cdot v$ the standard 
Euclidean hermitian inner product on $\C^n$, we obtain:
$$
\<y, y\>_h= \frac{(-1 +|z|^2) |y|^2-  (z\cdot y)(y\cdot z)}{-(-1 + |z|^2)^2}, \quad y\in T_z\B^n. 
$$
In the differential form, the metric $h$ is, thus, given by
$$
ds^2_h= \frac{1}{1-|z|^2} \sum_{k=1}^n dz_k d\bar{z}_k + \frac{1}{(1 - |z|^2)^2} \sum_{j, k=1}^n z_j \bar{z}_k dz_k d\bar{z}_j. 
$$
This hermitian metric is K\"ahler, with the K\"ahler potential (centered at the origin) equal to 
$$
f(z)= \log( 1- |z|^2), 
$$
and the K\"ahler form $\omega= \frac{i}{2} \D \bar\D f$ 
equal
$$
\omega = \frac{1}{1-|z|^2} \sum_{k=1}^n dz_k \wedge d\bar{z}_k + \frac{1}{(1 - |z|^2)^2} \sum_{j, k=1}^n z_j \bar{z}_k dz_k\wedge  d\bar{z}_j. 
$$

The complex hyperbolic metric on $\B^n$ (the unit ball in $\C^n$) is the Bergman metric with the Bergman kernel $K(z,\zeta)$ equal
$$
K(z,\zeta)= \frac{n!}{2\pi^n} (1- (z\cdot \zeta))^{-n-1},
$$
where, as before, $z\cdot \zeta$ is the standard hermitian inner product on $\C^n$. 
%Then the complex hyperbolic hermitian metric on $\B^n$ equals 
%$$h_z= \D \bar \D \log K(z,z). $$

The  distance function $d$ on $\H^n_\C$ satisfies 
$$
\cosh^2(d([x], [y]))= \frac{\<x,y\> \<x, y\>}{\<x,x\> \<y, y\>}. 
$$
For example, specializing to the case when $[x]$ is the center of $\B^n$ and $[y]$ is represented by a point $z\in \B^n$, we obtain:
$$
\cosh^2(d(0, z))= (1-|z|^2)^{-1}. 
$$

See \cite[pp. 72--79]{Gol} and \cite[\S 1.4]{Krantz}; note however that Goldman uses a different normalization of the metric on the complex hyperbolic space: with his normalization, the sectional curvature varies in the interval $[-2, -\frac{1}{2}]$. 
%(normalization?). 

\medskip 

A real linear subspace $W\subset V$ is said to be {\em totally real with respect to the form $\<\cdot, \cdot\>$} if for any two vectors 
$z, w\in W$,  $\<z, w\>\in \R$. Such a subspace is automatically totally real in the usual sense: $J W\cap W=\{0\}$, where $J$ is the almost complex structure on $V$. 

{\em Real geodesics} in $\B^n$ are projections (under $p$) of totally real indefinite (with respect to $q$) 2-planes in $V$ (intersected with $V_-$). For instance, geodesics through the origin $0\in \B^n$ are Euclidean line segments in $\B^n$.

More generally, totally-geodesic real subspaces in $\B^n$ are projections of totally real indefinite subspaces in $V$ (intersected with $V_-$). They are isometric to the real hyperbolic space $\H^n_\R$ of constant sectional curvature $-1$. Boundaries of real hyperbolic planes are called {\em real circles} in $S^{2n-1}$.

{\em Complex geodesics} in $\B^n$ are projections of indefinite complex 2-planes;  boundaries of complex geodesics are called {\em complex circles} in $S^{2n-1}$. 
Complex geodesics are isometric to the unit disk with the hermitian metric 
$$
\frac{dz d\bar z}{(1-|z|^2)^2},
$$ 
which has constant curvature $-4$. These are the extremal disks for the Kobayashi metric on $\B^n$, which coincides with the complex hyperbolic distance function $d$. It also equals the Caratheodory's metric on $\B^n$ (as is the case for all bounded convex domains in $\C^n$).

More generally, complex hyperbolic $k$-dimensional subspaces $\H^k_\C$ 
in $\B^n$ are projections of indefinite complex $k+1$-dimensional subspaces 
(intersected with $V_-$). 

All complete totally-geodesic submanifolds in $\H^n_\C$ are either real or complex hyperbolic subspaces.

The holomorphic bisectional curvature of $\H_\C^n$ is constant, equal $-1$. 
It turns out that $\H_\C^n$ has negative sectional curvature which varies in the interval $[-4,-1]$. Thus, $\H_\C^n$ is a {\em negatively pinched Hadamard manifold}:

\begin{defn}
1. A Hadamard manifold $X$ is a simply-connected complete nonpositively  curved Riemannian manifold. 

2. A Hadamard manifold $X$ is said to have {\em strictly negative curvature} if there exists $a< 0$ such that the sectional curvature of $X$ is $\le a$. 

3. A Hadamard manifold $X$ is said to be {\em  negatively pinched} (has {\em pinched negative curvature}) 
 if there exist two negative numbers $b\le a< 0$ such that the sectional curvature of $X$ lies in the interval $[b,a]$. 
\end{defn}

\medskip 
The group $U(n,1)=U(q)$ of (complex) automorphisms of $q$ projects to the group $G=PU(n,1)=Aut(\B^n)$ of complex (biholomorphic) automorphisms of $\B^n$. This group acts transitively, with the stabilizer of the center of $\B^n$ equal to $K=U(n)$. Consequently, the metric $d$ on $\B^n$ is complete. The group $G$ is a Lie group, its Lie topology is equivalent to the topology of pointwise convergence, equivalently, the topology of uniform convergence on compacts in $\B^n$, equivalently, the quotient topology of the matrix group topology on $U(n,1)$. The group $G$ is linear, its matrix representation is given, for instance, by the adjoint representation, which is faithful since $G$ has trivial center. 

The Lie group $G$ is connected and has real rank 1. Its  Cartan decomposition is 
$$
G= KA_+ K,
$$
where $A_+$ is the semigroup of positive translations (transvections) along a chosen geodesic through $0$. 

Let $\ol{\B^n}$ denote the closure of $\B^n$ in $\P^n$. The boundary sphere $S^{2n-1}= \partial \B^n$ of $\B^n$ is the projection to $\P^n$ of the null-cone of the form $q$. The sphere $S^{2n-1}$ a {\em CR manifold}: It is equipped with a smooth totally nonintegrable hyperplane distribution $H_z, z\in S^{2n-1}$,
$$
H_z= T_z S^{2n-1}\cap J(T_z S^{2n-1}),
$$
where $J$ is the almost complex structure on $\P^n$. The subspace $H_z$ is a (complex) hyperplane in $T_{z}\P^n$. We let $P_z$ denote the unique projective subspace in $\P^n$ passing through $z$ and tangent to $H_z$. Thus, $P_z\cap \ol{\B^n}=\{z\}$. 

One defines a {\em sub-Riemannian metric} $d_C$ on $S^{2n-1}$ as follows. Given points $\xi, \eta\in S^{2n-1}$, define $C_{p,q}$ as the collection of smooth paths $c: [0,1]\to S^{2n-1}$ connecting $p$ to $q$ such that $c$ is a {\em contact path}, i.e. $c'(t)\in H_{c(t)}$ for all $t\in [0,1]$. Then the {\em Carnot metric} $d_C$ on $S^{2n-1}$ is 
$$
d_C(\xi,\eta)=\inf_{c\in C_{\xi,\eta}} \int_0^1 ||c'(t)||dt,
$$
where $||\cdot||$ is a background Riemannian metric on $S^{2n-1}$, say, the unique metric of sectional curvature $+1$ invariant under 
$O(2n)$. It turns out that $d_C$ is indeed a metric which topologizes  $S^{2n-1}$.  However, unlike a Riemannian metric on $S^{2n-1}$, which has Hausdorff dimension equal to the topological dimension, 
the metric space $(S^{2n-1}, d_C)$ is {\em fractal}, its Hausdorff dimension
$\dim_H$ equals
$$
\dim_H (S^{2n-1}, d_C)=2n. 
$$ 

\medskip 
Most of the following discussion is valid for general negatively pinched Hadamard spaces; I refer to the paper by Bowditch \cite{Bowditch} for details, especially in the context of discrete isometry groups.

Since $\H^n_\C$ is a Hadamard manifold $X$, it has an {\em intrinsically defined} ideal (visual) boundary $\geo X$, defined as the set of equivalence classes of geodesic rays, where two rays are equivalent iff they are within finite Hausdorff distance. Every geodesic ray is equivalent to a geodesic ray emanating from a chosen base-point $o\in X$. The topology 
on $\geo X$ is the quotient topology, where the space of geodesic rays is equipped with the topology of uniform convergence on compacts. Equivalently, since the map from the unit tangent sphere $UT_oX$ at $o$ to $\geo X$ is bijective, $\geo X$ is homeomorphic to $UT_oX$. The union $\ol{X}:=X\cup \geo X$ also has a natural topology with respect to which it is homeomorphic to 
the closed ball.  Given a subset $Y\subset X$, we define $\geo Y$ as the intersection of the closure of $Y$ in $\ol{X}$ with $\geo X$. 

If  $X$ is strictly negatively curved, it satisfies the {\em visibility property}: Any two distinct points $\xi, \eta\in \geo X$ are connected by a unique geodesic, denoted $\xi\eta$.

In the case $X=\H_\C^n$, this abstract compactification is naturally homeomorphic to the closed ball compactification 
 $\ol{\B^n}$: Two geodesic rays $c_1, c_2$ are equivalent iff they terminate at the same point of the boundary sphere $S^{2n-1}$ (see \cite{Paupert} for details).

Suppose that $X$ is a Hadamard manifold. Given a closed subset $\La\subset \geo X$, one defines the {\em closed convex hull}, denoted $\hull(\La)$, of 
$\La$ in $X$ as the intersection of all closed subsets 
$C\subset X$ such that $\geo C\supset \La$. For $\eta>0$ we will use the notation $\hull_\eta(\La)$ to denote the closed $\eta$-neighborhood of $\hull(\La)$ in $X$. 

\begin{thm}
[M. Anderson, \cite{Anderson}] If $X$ has pinched negative curvature then for every closed subset $\La\subset \geo X$ which is not a singleton, $\hull(\La)$ is a (closed, convex) subset of $X$ such that 
$\geo \hull(\La)=\La$. 
\end{thm}

\begin{exe}
(a) Assuming that $X$ is negatively curved, verify: 

1. $\hull(\La)=\emptyset$ if and only if $\La$ consists of at most one point.

2. For any two distinct points $\xi, \eta\in \geo X$, $\hull(\{\xi, \eta\})= \xi \eta$. \

(b) Verify that Anderson's theorem fails for the Euclidean plane $X=E^2$. 
\end{exe} 

Anderson's theorem requires negative pinching: It fails if $X$ merely has strictly negative curvature, see \cite{Ancona}. 

The geometry of convex hulls remains a bit of a mystery, for instance we still do not entirely understand volumes of convex hulls of finite subsets. The best known result seems to be:

\begin{thm}
[A. Borb\'ely, \cite{Borbely}] If $X$ is $m$-dimensional, has curvature in the interval $[-k^2, -1]$  and $\La$ has cardinality $\le n$, then $Vol(\hull(\La))\le Cn^{2-\eta}$, where $C=C(m,k)$, while 
$$
\eta= \frac{1}{1+4k(m-1)}. 
$$
\end{thm}

\begin{defn}
For a closed subset $\La\subset \D \B^n$, define its {\em tangent hull} $\hat\La$ as the union of hyperplanes $P_\la, \la\in \La$. I will refer to the hyperplanes $P_\la, \la\in \La$ as the 
{\em complex support hyperplanes} of $\La$. Similarly, for an open subset $\Om= \D \B^n- \La$, define
\begin{equation}\label{eq:thull}
\check{\Om}= \P^n - \hat\La. 
\end{equation}
\end{defn}

\medskip 
\begin{exe}
$\hat\La$ is also closed and $\hat\La\cap \ol{\B^n}= \La$. 
\end{exe}

See Appendix A for a discussion of {\em horospheres} and {\em horoballs} in Hadamard manifolds $X$ and the {\em horofunction compactification} of $X$, which leads to an alternative description of the topology on $\ol{X}$. 

Isometries of $X$ extend to homeomorphisms of $\ol X$; in the setting of $\B^n$, this is just the fact that all automorphisms of $\B^n$ are restrictions of projective transformations:
$$
PU(n,1)< PGL(n+1, \C). 
$$

%There is a $G$-equivariant homeomorphism from $UT \B^n$ to the complement in the diagonal in $S^{2n-1}\times S^{2n-1}$: Given a unit tangent vector $v\in T_z\N^n$, we consider the unique 
%biinfinite geodesic $c$ in $\B^n$ such that $c(0)=z, c'(0)=v$. Then $c$ connects two distinct points $

The group $G=PU(n,1)$ acts doubly transitively on the boundary sphere $S^{2n-1}$: Given two pairs of distinct 
points $\xi_i, \eta_i, i=1,2$, we connect these points by unique biinfinite (unit speed) geodesics $c_i=\xi_i\eta_i$. Set $z_i:= c_i(0), v_i:=c'(0)\in T_{z_i}\B^n$. 
Then, since $G$ acts transitively on the unit tangent bundle $UT \B^n$, there exists $g\in G$ sending $v_1\mapsto v_2$. Thus, $g(c_1)=c_2$ and, consequently, 
$g(\xi_1)=\xi_2, g(\eta_1)=\eta_2$.

\medskip 
{\bf Classification of isometries.} Every isometry $g\in G=Aut(\B^n)$ is continuous on the closed ball $\ol{\B^n}$ and, hence, has at least one fixed point there. Accordingly, automorphisms $g\in G$ are classified as:

\begin{enumerate}
\item {\bf Elliptic}: $g$ has a fixed point $z$ in $\B^n$. After conjugating $g$ via $h\in Aut(\B^n)$ which sends $z$ to $0$,  
$$
hgh^{-1}\in K=U(n). 
$$
\item {\bf Parabolic}: $g$ has a unique fixed point in $\ol{\B^n}$ and this is a boundary point $z\in S^{2n-1}$. Equivalently, 
$$
\inf \{d(z, gz): z\in \B^n\}=0
$$
and the infimum is not realized. 

\item {\bf Hyperbolic}: $g$ has exactly two fixed points $\xi,\eta$  in $\ol{\B^n}$, both are in $S^{2n-1}$. (In particular, $g$ preserves the unique geodesic $\xi\eta$ in $\B^n$ and acts as a translation along this geodesic. This geodesic is called the {\em axis} of $g$.)  Equivalently, 
$$
\inf \{d(z, gz): z\in \B^n\}\ne 0. 
$$
This infimum is realized by any point on the axis of $g$. 
\end{enumerate}

The fixed point $\la$ of a hyperbolic isometry $\ga$ is called {\em attractive} (resp. {\em repulsive}) if for some (every) $x\in X$, $\ga^i(x)\to \la$ as $i\to\infty$ (resp. $i\to -\infty$).  

\medskip 
An elliptic automorphism of $\B^n$ is called a {\em complex reflection} if its fixed-point set is a complex hyperbolic hyperplane in $\H^n_\C$. 

 \medskip 
 As any strictly negatively curved Hadamard manifold, $\H^n_\C$ satisfies the {\em convergence property:}  
 
 \begin{thm}
 For every sequence $g_i\in G=PU(n,1)$, after extraction, the following dichotomy holds:
 
 (a) Either $g_i$ converges to an isometry $g\in G$. 
 
 (b) Or there is a pair of points $\xi, \eta\in S^{2n-1}$ such that $g_i|_{\ol{\B^n}-\{\eta\}}$ converges uniformly on compacts to the constant $\xi$. 
 \end{thm}
\proof ~First, prove this for sequences of hyperbolic isometries with a common axis. Then use the Cartan decomposition of $G$. \qed 

\medskip 
In the case (b), we will say that $(g_i)$ converges to the {\em quasiconstant map} $\xi_\eta$. (The point $\eta$ is the {\em indeterminacy point} of $\xi_\eta$.) 

It turns out that most elementary properties of discrete isometry groups of strictly negatively curved Hadamard manifolds can be derived just from the Convergence Property! 
See \cite{Bowditch1999, Tukia1994, Tukia1998} for a development of the theory of 
{\em convergence group actions} on compact metrizable spaces, i.e. topological group actions satisfying the Convergence Property. 

\begin{exe}
1. Verify that if $g_i\to \xi_\eta$ then $g_i^{-1}\to \eta_\xi$. 

2. If $g_i\to \xi_\eta$, verify that $(g_i)$ converges (again, uniformly on compacts) to the constant map $\xi$ on $\P^n - P_\eta$. 

3. Find an example where $\xi=\eta$.
\end{exe}

%\section{Further structure}

%Symplectic structure. Symplectomorphisms. Quasiisometries. Contactomorphisms and quasiconformal maps of the boundary sphere. 

\section{Basics of discrete subgroups of $PU(n,1)$}

Almost all the properties of discrete subgroups $\Ga< G=PU(n,1)$ stated in this section hold for discrete isometry groups of negatively pinched Hadamard manifolds. 

\begin{defn}
A subgroup $\Ga< \Isom(X)$ of isometries of a Riemannian manifold $X$ is called {\em discrete} if it is discrete as a subset of $\Isom(X)$.  
Discrete subgroups $\Ga< PU(n,1)$ are {\em complex hyperbolic Kleinian groups}. 
%A complex hyperbolic manifold/orbifold is the quotient manifold/orbifold $M_\Ga=\H^n_\C/\Ga$ of a complex hyperbolic Kleinian group $\Ga< PU(n,1)$.  
\end{defn}

Here, all {\em reasonable} topologies on $\Isom(X)$ agree. For instance, one can use the topology of uniform convergence on compact subsets, or the topology of pointwise convergence. 
%In the case of $X=\H^n_\C$, the isometry group of $X$ is linear and the above topology is equivalent to the standard matrix topology. 

Recall that a group $\Ga$ of homeomorphisms of a Hausdorff 
topological space $X$ is said to  act {\em properly discontinuously} 
on $X$ if for every compact $C\subset X$,
$$
\card \{\ga\in \Ga: \ga C\cap C\ne \emptyset\} <\infty. 
$$
See \cite{proper} for a comparison of alternative notions of proper discontinuity.

\begin{exe}
Suppose that $X$ is a Riemannian manifold and $G=\Isom(X)$ is the isometry group of $X$. 

(a) Prove that the following are equivalent for subgroups $\Ga< G$:

1. $\Ga$ is a discrete subgroup  of $G$. 

2. $\Ga$ acts properly discontinuously on $X$. 

3. For one (equivalently, every) $x\in X$ the function $\Ga\to \R_+$, $\ga\mapsto d(x, \ga x)$ is proper (with $\Ga$ equipped with discrete topology), 
i.e. if $\ga_i$ is a sequence consisting of distinct elements of $\Ga$, then
$$
\lim_{i\to\infty} d(x, \ga_i x)=\infty. 
$$

(b) Every discrete subgroup of $G$ is at most countable. 
\end{exe}

A group $\Ga$ is said to act {\em freely} on $X$ is for every $x\in X$, the $\Ga$-stabilizer 
$$\Ga_x=\{\ga\in \Ga: \ga x=x\}$$ is the trivial subgroup of $\Ga$. 

If $X$ is a manifold and $\Ga$ is a group acting freely and properly discontinuously, then the quotient space $X/\Ga$ is a manifold and the projection map $X\to X/\Ga$ is a covering map. 
If one does not assume freeness of the action then $X/\Ga$ is an {\em orbifold} and the projection map $X\to X/\Ga$ is an {\em orbi-covering map}.  If $X$ is simply-connected, the group $\Ga$ is the (orbifold) fundamental group of $X/\Ga$. See Appendix D for a discussion of orbifolds and related concepts.

In the case when $X$ is a Hadamard manifold, 
a  subgroup $\Ga< \Isom(X)$ acts freely on $X$ if and only if $\Ga$ is {\em torsion-free}, i.e. every nontrivial element of $\Ga$ has infinite order. If $\Ga$ acts on $X$ isometrically/holomorphically, 
the Riemannian metric/complex structure on $X$ descends to the quotient manifold (orbifold) $X/\Ga$.

\begin{defn}
A complex hyperbolic $n$-dimensional orbifold (manifold) is the quotient of  $\H^n_\C$ by a discrete (torsion-free) subgroup of $PU(n,1)$, $M_\Ga=\H^n_\C/\Ga$.  
\end{defn}

\begin{exe}
Assuming that $X$ is a Hadamard manifold and $\Ga<\Isom(X)$ is discrete, prove that $\Ga$ is torsion-free if and only if it contains no elliptic elements, besides the identity. 
\end{exe}

For {\em finitely generated subgroups} $\Ga< PU(n,1)$, one can eliminate torsion by passing to a finite index subgroup:

\begin{thm}
[Selberg's Lemma, see e.g. \cite{DK} or \cite{Ratcliffe}] If $\k$ is a field and $\Ga< GL(n, \k)$ is a finitely generated subgroup, then $\Ga$ is {\em virtually torsion-free}, i.e. contains a torsion-free subgroup of finite index.  
\end{thm}

In particular, every complex hyperbolic orbifold $\O$ with finitely generated (orbifold) fundamental group, admits a finite-sheeted manifold orbi-covering $M\to \O$.

\begin{rem}
Selberg's theorem fails for discrete finitely generated groups of isometries of negatively pinched Hadamard manifolds, see \cite{K2018}. 
\end{rem}

\begin{defn}
Given a Hadamard manifold $X$, a discrete subgroup $\Ga< \Isom(X)$ and a point $x\in X$, the {\em limit set} $\La=\La_{\Ga}$ is the accumulation set of the orbit $\Ga x$ in $\geo X$, i.e. 
$$
\La= \geo (\Ga x). 
$$ 
The complement $\Om:= \geo X - \La$ is called the {\em discontinuity domain} of $\Ga$. 
\end{defn}

\begin{exe}\label{exe:limit}
Suppose that $\Ga$ is a discrete subgroup of $\Isom(X)$ and $X$ is strictly negatively curved. Verify:

\begin{enumerate}
\item $\La$ is independent of $x\in X$. (Hint: Use the Convergence Property.)\footnote{This also holds for general Hadamard manifolds  even though the convergence property fails.}  

\item $\La$ is closed and $\Ga$-invariant.  Accordingly, $\Om$ is open in $\geo X$ and is $\Ga$-invariant as well. 

%\item Either $\La=\geo X$ or $\La$ has empty interior. 

\item $\Om$ is either empty or is dense in $\geo X$.  (Hint: Use the Convergence Property.)

\item Either $\La$ consists of at most two points or it is perfect, i.e. contains no isolated points.  

\item If $\Ga'$ is a subgroup of $\Ga$, then $\La_{\Ga'}\subset \La_{\Ga}$. 

\item If $\Ga'\triangleleft \Ga$ is an infinite normal subgroup then $\La_{\Ga'}= \La_{\Ga}$.

\item If $\Ga'< \Ga$ is a subgroup of finite index then $\La_{\Ga'}= \La_{\Ga}$. 
\end{enumerate}
\end{exe}

\begin{example}
Let $\ga\in \Isom(X)$ be a non-elliptic element. Then the limit set of the group $\Ga=\<\ga\>$ generated by $\ga$ equals the fixed-point set of $\ga$ in $\geo X$. 
\end{example}

\begin{lemma}
If $\Ga< \Isom(X)$ is a discrete subgroup and $X$ is a strictly negatively curved Hadamard manifold, then 
$\Ga$ acts properly discontinuously on $Y=X\cup \Om$. 
\end{lemma}
\proof ~Let $C$ be a compact subset of $Y$. Suppose that there exists a sequence consisting of distinct elements $\ga_i\in \Ga$ such that for each $i$, $\ga_i C\cap C\ne \emptyset$. In view of the Convergence Property, after extraction, the sequence $\ga_i$ either converges to an isometry $\ga\in \Isom(X)$ (which would contradict the discreteness of $\Ga$) or to a quasiconstant map $\xi_\eta$, with $\xi, \eta\in \La$. Since $(\ga_i)$ converges to $\xi$ uniformly on compacts in $\ol{X}- \{\eta\}$ and $C\subset Y \subset \ol{X}- \{\eta\}$ is compact, there exists a neighborhood $U$ of $\xi$ disjoint from $C$; thus, for all but finitely many values of $i$, $\ga_i(C)\subset U$. A contradiction. \qed

\medskip 
A more difficult result is

\begin{thm}
[A. Cano, J. Seade,  \cite{CS, CNS}]  Every discrete subgroup $\Ga< PU(n,1)$ acts properly discontinuously on $\check{\Om}:= \P^n - \hat\La$ (see \eqref{eq:thull}).  
\end{thm}

\begin{rem}\label{rem:pd} 
An alternative proof of this result is an application of a more general proper discontinuity theorem proven in \cite{KLP}. More precisely, let $F_{1,n}$ be the flag-manifold consisting of flags $(V_1, V_n)$ in $V=\C^{n+1}$, where $V_1$ is a line and $V_n$ is a hyperplane (containing $V_1$). We have a $G$-equivariant holomorphic fibration $\pi: F_{1,n}\to \P^n$ sending each pair $(V_1, V_n)$ to $V_1$. 
The tangent hull $\hat\La$ of $\La$ defines a natural continuous map $\theta: \La\to F_{1,n}$ sending each $\la\in \La$ to the pair $(V_1, V_n)$ consisting of the preimages of $\la$ and $P_\la$ in $V$. Let $\tilde\La$ be the image of $\theta$ and let $Th(\tilde\La)$ be the {\em thickening} of  $\tilde\La$ in $F_{1,n}$, consisting of flags 
$(V_1', V_n')$ such that either  $V_1'$ belongs to $\La$ or $V_n'$ is a complex support hyperplane of $\La$. Then $\Ga$ acts properly discontinuously on $\Om_{Th}= F_{1,n} - Th(\tilde\La)$; see 
\cite{KLP}. Since $\pi^{-1}(\check{\Om})\subset \Om_{Th}$, the action of $\Ga$ on  $\check{\Om}$ is properly discontinuous as well. 
\end{rem}

\medskip 
In particular, the quotient $\ol{M}_\Ga:=(\B^n \cup\Om)/\Ga$ embeds as an orbifold with boundary in the complex orbifold without boundary $\check{\Om}/\Ga$. The boundary of 
$\ol{M}_\Ga$ (equal to $\Om/\Ga$) is strictly Levi-convex in   $\check{\Om}/\Ga$. 

\begin{notation}
The {\em boundary} $\D M_\Ga$ of a complex hyperbolic orbifold $M_\Ga$ is $\Om_\Ga/\Ga$; in other words, this is the boundary of $\ol{M}_\Ga$. 
\end{notation}

We now return to the discussion of discrete subgroups of general negatively pinched Hadamard manifolds $X$. 

\begin{thm}
If $\al, \beta$ are hyperbolic elements of a discrete subgroup of $\Isom(X)$, then their fixed-point sets are either equal or disjoint. 
\end{thm}

\begin{cor}
If $\Ga< \Isom(X)$ is discrete and fixes a point $\la\in \geo X$ then $\La_{\Ga}$ either equals to $\{\la\}$ and $\Ga$ contains no hyperbolic elements, or $\La_{\Ga}$ consists of two points, $\La_\Ga=\{\la, \la'\}$ and $\Ga$ contains no parabolic elements.   
\end{cor}

\begin{defn}
A discrete subgroup $\Ga< \Isom(X)$ is called {\em elementary} if 
$\card (\La_\Ga)\le 2$. It is said to be {\em nonelementary} otherwise. 
\end{defn}

Elementary subgroups are, in many ways, exceptional, among discrete subgroups. 

\medskip 
In view of Exercise \ref{exe:limit}(4),  the limit set of every nonelementary subgroup is perfect. In particular, it has the cardinality of continuum. Hence:

\begin{prop}
The limit set of a discrete subgroup of $\Isom(X)$ consists of 0, 1, 2 or continuum of points. 
\end{prop}

\begin{prop}
The limit set of a nonelementary discrete group $\Ga$ is the smallest nonempty closed $\Ga$-invariant subset of $\geo X$.  In particular, every orbit  in $\La_\Ga$ is dense.  
\end{prop}
\proof Suppose that $L\subsetneq \La_\Ga$ is a closed nonempty and $\Ga$-invariant subset.  Take a point $\xi\in \La_\Ga -L$ and let $(\ga_i)$ be a sequence in $\Ga$ converging to a quasiconstant map $\xi_\eta$. Then, for every $\la\in L-\{\eta\}$, $\lim_{i\to\infty} \ga_i(\la)= \xi$. Since $L$ is closed and $\xi\notin L$, for all sufficiently large $i$, $\ga_i(\la)\notin L$, contradicting invariance of $L$. This leaves us with the possibility that $L$ is the singleton $\{\xi\}$  fixed by the entire $\Ga$. Then $\Ga$ is elementary.    
\qed

\begin{thm}\label{thm:elementary} 
Suppose that $\Ga$ is an elementary subgroup of $\Isom(X)$. 

1. If $\La_\Ga$ is a singleton then every element of $\Ga$ is elliptic or parabolic.

2.  If $\La_\Ga$ consists of two points then every element of $\Ga$ is elliptic or hyperbolic. Hyperbolic elements fix $\La_\Ga$ pointwise. Elliptic elements can swap the two limit points. 

3. $\Ga$ is a virtually nilpotent\footnote{I.e. contains a nilpotent subgroup of finite index} group. 
\end{thm}

See \cite{BK} for a more detailed discussion of elementary groups and their quotient spaces $M_\Ga$. Here we only note that 
discrete elementary subgroups of $PU(n,1)$ are virtually {\em 2-step}  nilpotent.

\begin{prop}
Suppose that $X$ is a strictly negatively curved Hadamard manifold. 
If $\xi, \eta$ are distinct limit point of a discrete subgroup $\Ga< \Isom(X)$, then there exists a sequence $\ga_k\in \Ga$ of hyperbolic elements whose attractive (resp. repulsive) fixed points converge to $\xi$ (resp. $\eta$). 
\end{prop}
\proof Since $\xi, \eta$ are limit points of $\Ga$, there exist sequences $(g_i), (h_j)$ in $\Ga$ which converge, respectively, to the quasiconstant maps $\xi_\al$ and $\beta_\eta$. By precomposing these sequences with  suitable 
elements of $\Ga$, we can assume that the points $\xi, \eta, \al, \beta$ are pairwise distinct. Let $U_\al, U_\beta, U_\xi, U_\eta$ be pairwise disjoint open ball neighborhoods in $\geo X$ 
of $\al, \beta, \xi, \eta$ respectively. In view of the convergence $g_i\to \xi_\al, h_j\to \beta_\eta$, for all sufficiently large $i$  we have 
$$
h_i(\geo X - U_\eta)\subset U_\beta, \quad g_i(\geo X  - U_\al)\subset U_\xi,$$
and, hence,
$$
g_i \circ h_i (\geo X - U_\eta)\subset U_\xi. 
$$ 
In particular, the composition $f_i=g_i \circ h_i$ has an attractive fixed point in $U_\xi$. Similarly, $f_i^{-1}$ has an attractive fixed point in $U_\eta$. 
\qed

\begin{cor}
If $\Ga$ is nonelementary then the set of hyperbolic fixed points of elements of $\Ga$ is dense in $\La_\Ga$. 
\end{cor}

\begin{cor}
If a discrete group $\Ga$ contains a parabolic element, then parabolic fixed points are dense in $\La_\Ga$.  
\end{cor}

The following theorem provides a converse to Theorem \ref{thm:elementary}(3): 

\begin{thm}
Each nonelementary discrete subgroup $\Ga< \Isom(X)$ contains a nonabelian free subgroup whose limit set is homeomorphic to the Cantor set. 
\end{thm}

\begin{defn}
The {\em convex core}, $Core(M)$, of $M=M_\Ga=X/\Ga$ is the projection to $M_\Ga$ of the closed convex hull $\hull(\La_\Ga)$ 
of the limit set of $\Ga$. 
\end{defn}

Given $\eta>0$, define  $Core_\eta(M)$ as the projection to $M_\Ga$ of $\hull_\eta(\La_\Ga)$. Intrinsically, the convex core can be defined as: 

\begin{exe}
$Core(M)$ is the intersection of all closed convex suborbifolds $M'\subset M$ such that $\pi_1(M')\to \pi_1(M)$ is surjective. 
\end{exe}

\medskip 
{\bf Conical limit points.} I conclude this section with a discussion of 
a  classification of limit points of discrete subgroups of $\Isom(X)$. 

\begin{defn}
A sequence $(x_i)$ in $X$ is said to converge to a point $\xi\in \geo X$ {\em conically} if there exists a geodesic ray $x\xi$ in $X$ and a constant $R<\infty$ such that:

$d(x_i, x\xi)\le R$ for all $i$ and $\lim_{i\to\infty} x_i=\xi$. 
\end{defn}

\begin{exe}
Let $\la\in \La_\Ga$ be a limit point. %Consider a geodesic ray $x\la\subset X$. 
The following are equivalent:

1. There exists a  sequence $\ga_i\in \Ga$ such that the sequence $(\ga_i(x))$ converges to $\xi$ conically. 

2. The projection of the ray $x\la$ to $M_\Ga$ defines a non-proper map $\R_+\to M_\Ga$. 
\end{exe}

\begin{defn}
A limit point $\la\in \La_\Ga$ is called {\em conical} or {\em radial} if it satisfies one of the two equivalent properties in this exercise. The set of conical limit points of $\Ga$ is denoted $\La^c=\La^c_\Ga$. 
\end{defn}

\begin{example}
1. If $\Ga$ is an elementary hyperbolic subgroup of $\Isom(X)$ then $\La_\Ga= \La^c_\Ga$.

2. If  $\Ga$ is an elementary parabolic subgroup of $\Isom(X)$ then $\La^c_\Ga=\emptyset$. 
\end{example}

\section{Margulis Lemma and thick-thin decomposition}

In this section, $X$ is a negatively pinched Hadamard manifold. 
For each discrete subgroup $\Ga< \Isom(X)$, a point $x\in X$ and a number $\eps>0$, define $\Ga_{x,\eps}$ to be the subgroup of $\Ga$ generated by the (necessarily finite) set
$$
\{\ga\in \Ga: d(x, \ga x)<\eps\}. 
$$
This subgroup is the ``almost-stabilizer'' of $x$ in $\Ga$. 

Let $U_{\Ga,\eps}$ denote the subset of $X$ consisting of points $x$ for which the almost-stabilizer  $\Ga_{x,\eps}$  is infinite. 

The components of $U_{\Ga,\eps}$ need not be convex (already for $X=\H^2_\C$), but each component is contractible: 

\begin{prop}
Each component of  $U_{\Ga,\eps}$ is contractible. 
\end{prop}

 In view of the contractibility of $X$ and of $\hull \La_\Ga$, it follows that $X- U_{\Ga,\eps}$ and 
$\hull \La_\Ga - U_{\Ga,\eps}$ are both contractible. Furthermore, if $X$ has curvature $\le -1$, each component $U$ of $U_{\Ga,\eps}$ is {\em uniformly quasiconvex}:

\begin{thm}
There exist universal constants $\delta_0, \eta_0$ such that 
each component $U$ of $U_{\Ga,\eps}$ satisfies:

1. For any two points $x, y\in U$, the geodesic $xy$ is contained in the $\delta_0$-neighborhood of $U$. 

2.  The $\eta_0$-neighborhood of $U$ is convex.  
\end{thm}

\begin{thm}[Kazhdan--Margulis; Margulis; see e.g. \cite{BGS}] There exists $\eps=\eps(n,b)$ such that every 
$n$-dimensional Hadamard manifold $X$ of sectional curvature bounded below by $b\le 0$ satisfies the following. 
For every discrete subgroup $\Ga< \Isom(X)$ and every $x\in X$, the subgroup $\Ga_{x,\eps}$ is virtually nilpotent. In particular, if $X$ is negatively curved, then  $\Ga_{x,\eps}$ is elementary. 
\end{thm}

\begin{cor}
For each discrete torsion-free subgroup $\Ga< \Isom(X)$,  the set $U_{\Ga,\eps}$ breaks into connected components $X_{\Ga,\eps,i}$ each of which is stabilized by some elementary subgroup $\Ga_i$ 
of $\Ga$ and for each $x\in X_{\Ga,\eps,i}$ the stabilizer $\Ga_i$ contains the 
``almost stabilizer''  $\Ga_{x,\eps}$. (The index can be infinite.)  
\end{cor}

As a corollary, one obtains the {\em thick-thin decomposition} of the orbifold $M=M_\Ga$:

\medskip 
\noindent $M_{(0,\eps)}$ is the projection of  $U_{\Ga,\eps}$ to $M$. It consists of all points $y\in M$ for which there exists a homotopically nontrivial  loop based at $y$ of length $< \eps$. Define also $M_{(0,\eps]}$ as the closure of $M_{(0,\eps)}$ in $M$. Both 
$M_{(0,\eps)}$ and $M_{(0,\eps]}$ are called  the $\eps$-{\em thin} parts of $M$. 
The complement $M_{[\eps,\infty)}= M - M_{(0,\eps)}$ and its interior $M_{(\eps,\infty)}$ 
are called the {\em $\eps$-thick} parts of $M$. 

One defines the $\eps$-thick, resp. thin, part of the convex core $Core(M)$ as the intersection $Core(M)\cap M_{[\eps,\infty)}$, resp. $Core(M)\cap M_{(0,\eps)}$. 

\medskip 
Components of the thin parts $M$ and $Core(M)$ come in two shapes: 

(a) {\bf Tubes}.  Suppose that $U$ is a component of $U_{\Ga,\eps}$ whose stabilizer $\Ga_U$ in $\Ga$ is virtually hyperbolic, i.e. contains a cyclic hyperbolic subgroup of finite index.  In other words, the limit set of 
$\Ga_U$ consists of two points $\xi, \eta$. The geodesic $\xi\eta$ is then invariant under $\Ga_U$; it is also contained in $U$ and projects to a closed geodesic $c\subset U/\Ga_U$. 
The quotient $U/\Ga_U$ is a {\em tube}: If $\Ga_U$ is torsion-free then this quotient is homeomorphic to an $\R^k$-bundle over $S^1$, with the base of the fibration corresponding  
to  the closed geodesic $c$.

(b) {\bf Cusps}.  Suppose that $U$ is a component of $U_{\Ga,\eps}$ whose stabilizer $\Ga_U$ in $\Ga$ is virtually parabolic, i.e. contains a parabolic subgroup of finite index. 
 In other words, the limit set of $\Ga_U$ consists of a single point $\eta$. The group $\Ga_U$ preserves  horoballs $B_\eta$ based at $\eta$. The subsets  $U_{\Ga,\eps}$ are typically strictly smaller (not even Hausdorff-close) than any of the horoballs $B_\eta$.

\section{Geometrically finite groups}

The notion of geometrically finite Kleinian groups was introduced by Lars Ahlfors in mid 1960s for the real hyperbolic space 
and later generalized (by William Thurston and Brian Bowditch) to manifolds of negative curvature: The discrete groups in this class are the {\em nicest-behaving} among discrete isometry groups of negatively pinched Hadamard manifolds. 

\begin{defn}
Let $X$ be a negatively pinched Hadamard manifold. 
A discrete subgroup $\Ga< G=\Isom(X)$ is called {\em geometrically finite} if:

(a) The orders of elliptic elements of $\Ga$ are uniformly bounded (from above), and 

(b) the volume of  $Core_\eta(M_\Ga)$ is finite for some (equivalently, every, $\eta>0$).  

\noindent A  discrete subgroup $\Ga< G$ is called {\em convex-cocompact} if $\card (\La_\Ga)\ne 1$ and $Core(M_\Ga)$ is compact. 
\end{defn}

For instance, if $\La_\Ga=\geo X$ then $\hull(\La_\Ga)=X$ and, thus,  $\Ga$ is geometrically finite if and only if $\Ga< G$ is a {\em lattice}, i.e. $vol(M_\Ga)<\infty$. 
Under the same assumption, $\Ga$ is convex-cocompact if and only if $\Ga< G$ is a {\em uniform lattice}, i.e. $M_\Ga$ is compact. 

\begin{thm}
1. (B. Bowditch, \cite{Bowditch}) A discrete subgroup $\Ga< G$ is geometrically finite iff  the $\eps$-thick part of $Core(M_\Ga)$ is compact. 

2. (B. Bowditch, \cite{Bowditch}) A  discrete subgroup $\Ga< G$ is convex-cocompact iff $\ol{M}_\Ga$ is compact. 

3. (B. Bowditch, \cite{Bowditch}) A  discrete subgroup $\Ga< G$ is convex-cocompact iff every limit point of $\Ga$ is conical. 

4. (M. Kapovich, B. Liu, \cite{KL}) A discrete subgroup $\Ga< G$ is geometrically finite iff every limit point of $\Ga$ is either conical or a parabolic fixed point.  
\end{thm}

In particular, (1) implies that geometrically finite groups are finitely presentable (since $\hull \La_\Ga - U_{\Ga,\eps}$ is contractible).

In particular, a convex-cocompact subgroup $\Ga< PU(n,1)$ acts properly discontinuously an cocompactly on $\H^n_\C \cup \Om$. The action of $\Ga$ on 
 $\check{\Om}$ is properly discontinuous but (in general) not cocompact. If becomes cocompact if we lift to the flag-manifold $F_{1,n}$ (see \cite{KLP}): 

\begin{thm}
The $\Ga$-action on the domain $\Om_{Th}\subset F_{1,n}$ is properly discontinuous and cocompact. 
\end{thm}

\section{Ends of negatively curved manifolds}

Let $X$ be a negatively pinched Hadamard manifold and let 
$\La$ be a closed subset of $\geo X$ consisting of at least two points. Set $\Om= \geo X -\La$. 
The nearest-point projection $\Pi: X\to \hull(\La)$ extends continuously to a map $\Pi: X\cup \Om\to \hull(\La)$: 
While for $x\in X$, $\Pi(x)$ 
is defined by minimizing the distance function $d_x=d(x, \cdot)$ on $\hull(\La)$, for $\xi\in \Om$, the projection $\Pi(\xi)$ 
is defined by minimizing the Busemann function $b_\xi$ based at $\xi$. For a component $\Om_0\subset \Om$ we define a subset 
$X_0\subset X$  as the union of open %\footnote{i.e. with the point $x$ excluded} 
geodesic rays $x\xi -\{x\}$, where $\xi\in \Om_0, x=\Pi(\xi)$. The union of these geodesic rays is an open subset of $X-\hull(\La)$ whose closure in $X\cup\Om$ equals $X_0\cup \Om_0\cup \Pi(\Om_0)$. 

We now specialize to the setting when $\La=\La_\Ga$ is the limit set of a discrete subgroup $\Ga< \Isom(X)$. 
If $\Om_0$ has cocompact stabilizer $\Ga_0$ in $\Ga$, then $\Ga_0$ also acts cocompactly on  $X_0\cup \Om_0\cup \Pi(\Om_0)$. Thus, $M_\Ga$ has an the {\em isolated end} $E_0$ corresponding to $\Om_0/\Ga_0$,  
with the isolating neighborhood $X_0/\Ga_0$.  

\begin{defn}
Ends $E_0$ of $M=M_\Ga$ which have this form are called {\em convex ends} of $M$. 
\end{defn}

From the analytical viewpoint, the advantage of working with convex ends $E_0$ is that they admit {\em convex exhaustion functions}: For every convex end $E_0$ there exists a convex function $\phi: M\to \R_+$ which is proper on the closure of $E_0$ and vanishes on $M- E_0$. This function is the projection of the distance function $d(x, \hull_\eta(\La))$. 

\medskip 
Suppose that $C$ is an unbounded 
component  of the thin part $M_{(0,\eps)}$ of $M=M_\Ga$, and  $C$ has compact boundary. Then $C$ also defines 
an isolated end $E_C$ with an isolating neighborhood given by $C\cap M_{(0,\eps)}$.

\begin{defn}
Ends $E_C$ of $M_\Ga$ which have this form are called {\em cuspidal ends} of $M_\Ga$. 
\end{defn}

\begin{exe}
1. $\Ga$ is convex-cocompact iff $M_\Ga$ has only convex ends. 

2. If $M_\Ga$ has only convex and cuspidal ends then $\Ga$ is 
geometrically finite. 
\end{exe}

One can refine (cf. \cite{K2000}) the above definitions in two ways:

(a) Considering unbounded components of the thin part of $Core(M_\Ga)$ and, thus, defining cuspidal ends of the convex core. 

(b) Removing from $M_\Ga$ its cuspidal ends and their preimages under the nearest-point projection $M_\Ga\to Core(M_\Ga)$, one then defines {\em relative} convex ends of $M_\Ga$. 

\medskip 
One can also classify the ends of $M_\Ga$ using the potential theory as {\em hyperbolic} and {\em parabolic} ends, see \cite{NR1}. Note that if $M=M_\Ga$ is a complex hyperbolic manifold, then every convex end $E$ of $M$ is hyperbolic. 

\section{Critical exponent}

\begin{notation}
Let $B(x,r)$ denote the open ball of radius $r$ and center at $x$ in a metric space. 
\end{notation} 

The {\em critical exponent} of a discrete isometry group $\Ga$ of a Hadamard manifold $X$ (typically, satisfying some further curvature restrictions) 
is, probably, the single most important numerical invariant  of $\Ga$: It reflects both geometry of $\Ga$-orbits in $X$, the geometry of the limit set of $\Ga$, the ergodic theory of 
the action of $\Ga$ on the limit set and  analytic properties of the quotient space $X/\Ga$. Its origin goes back to the 19th century and the work of Poincar\'e (among others), who was interested in constructing {\em automorphic functions} (and forms) 
on the hyperbolic plane by ``averaging'' a certain holomorphic function (or a form) over a discrete isometry group $\Ga$. The resulting infinite series (the {\em Poincar\'e series}) may or may not converge, depending on the {\em weight} of the form, leading to the notion of the {\em critical exponent} or the {\em exponent of convergence} of $\Ga$. 

\medskip 
Let $\Ga< \Isom(X)$, a discrete isometry group of a Hadamard manifold. Pick points $x, y\in X$. The  {\em entropy} of $\Ga$ is defined as
$$
\delta=\delta_\Ga= \lim \sup_{r\to \infty} \frac{1}{r} \card (B(x, r) \cap \Ga y). 
$$
Thus, the entropy 
%critical exponent 
measures the rate of exponential growth of $\Ga$-orbits in $X$. It turns out that  $\delta$ equals the 
{\em critical exponent} of $\Ga$, defined as  
$$
\delta=\inf \{s: \sum_{\ga\in \Ga} \exp( -s d(x, \ga y)) < \infty\},
$$
i.e. $\delta$ is the {\em exponent of convergence} of the Poincare series $\sum_{\ga\in \Ga} \exp( -s d(x, \ga y))$.  Furthermore, 
$\delta$ is independent of the choice of $x, y\in X$. If 
$$
\sum_{\ga\in \Ga} \exp( -\delta d(x, \ga y)) < \infty
$$
(which depends only on $\Ga$ and not on the choice of $x, y$), then $\Ga$ is said to be a group of {\em convergence type}; otherwise, $\Ga$  is said to be of {\em divergence type}. 

\medskip 
Below are equivalent characterizations of $\delta$ in the case $X=\H^n_\C$: % both are (mostly) due to Corlette \cite{Corlette}:

\begin{thm}
Suppose that $\Ga< PU(n,1)$ is a discrete subgroup.Then: 

1. (Corlette \cite{Corlette}; Corlette--Iozzi \cite{CO}, Theorem 6.1) %If $\Ga$ is convex-cocompact %needed, or just conical limit set ? then 
$\delta=\delta_\Ga$ equals the Hausdorff dimension $\dim_H \La_\Ga^c$, where the conical limit set $\La^c_\Ga$ is equipped with the restriction of the Carnot metric on $S^{2n-1}$. In particular, if $\Ga$ is geometrically finite then $\delta= \dim_H \La$. 

2. (Elstrodt--Patterson--Sullivan--Corlette--Leuzinger, see \cite[Corollary 1]{Leuz}) Let $\la=\la(M_\Ga)$ denote the bottom of the $L^2$-spectrum of the Laplacian on $M_\Ga$. Then 
$$
\begin{cases}
\la= n^2 & \hbox{if~} 0\le \delta\le n\\ 
\la= \delta(2n-\delta) & \hbox{if~} n\le \delta\le 2n 
\end{cases} 
$$ 
\end{thm}

For further discussion of the critical exponent in the broader context of negatively curved Hadamard manifolds and Gromov-hyperbolic spaces, see e.g. \cite{Coo, DOP, DSU, Ledrappier, RT}.

\section{Examples}\label{sec:example} 

In this section I collect some examples of discrete subgroups of $PU(n,1)$ and their quotient-orbifolds. 

\begin{notation}
I will say that a discrete torsion-free subgroup $\Ga< G=PU(n,1)$ is {\em Stein} if the complex manifold $M_\Ga$ is Stein. 
\end{notation}

\medskip 
I will start with two elementary examples. 

\begin{example}
{\bf Cyclic hyperbolic groups.} Let $\ga\in PU(n,1)$ be a hyperbolic isometry fixing points $\la_\pm\in S^{2n-1}=\geo \H^n_\C$ 
and let $\Ga=\<\ga\>$ be the cyclic subgroup of $PU(n,1)$  generated by $\gamma$. Then  $\Ga< G$ is an elementary subgroup with the limit set $\La=\{\la_-, \la_+\}$. The quotient manifold 
$M_\Ga=\H^n_\C/\Ga$ is diffeomorphic to the product $\R^{2n-1}\times S^1$ while  $\ol{M}_\Ga$ is diffeomorphic to the product $\bar{D}^{2n-1}\times S^1$, where 
$\bar{D}^{2n-1}$ is the closed disk of real dimension $2n-1$. 
\end{example}

\begin{example}
{\bf Integer Heisenberg groups.} Given a natural number $n$, define the $2n+1$-dimensional real Lie group $H_{2n+1}$ as the group of $(n+2)\times (n+2)$-matrices
$$
\left[\begin{array}{ccc}
1 & \mathbf{a} & c\\ 
0 & I_{n} & \mathbf{b}\\
0 & 0 & 1
\end{array}\right],
$$
where $I_n$ is the identity $n\times n$ matrix, $\mathbf{a}\in \R^n$ is a row-vector, $\mathbf{b}\in \R^n$ is a column-vector and $c\in \R$. This group is 2-step nilpotent with the 1-dimensional center consisting of the matrices with 
$\mathbf{a}= \mathbf{b}=0$ and $c\in \R$. The quotient of $H_{2n+1}$ by its center is the $2n$-dimensional commutative Lie group isomorphic to $\R^{2n}$. The real Heisenberg group 
$H_{2n+1}$ contains the integer Heisenberg group $H_{2n+1}(\Z)$, defined as the intersection
$$
H_{2n+1}\cap SL(n+2, \Z). 
$$ 
The quotient $N=H_{2n+1}/H_{2n+1}(\Z)$ is a compact {\em nil-manifold}, which is a nontrivial circle over the torus $T^{2n}$. Algebraically, in terms of its presentation, 
$H_{2n+1}(\Z)$ is given by
$$
\< x_1, y_1,...,x_n, y_n, t| [x_i, t]=[y_j, t]=1, [x_i, y_i]=t, i=1,...,n, j=1,...,n\>.
$$

  The  Heisenberg group $H_{2n+1}$ embeds in $PU(n+1,1)$, fixing a point $\xi$ in $\geo \H^{n+1}_\C$ and acting simply-transitively on every horosphere in $\H^{n+1}_\C$ centered at $\xi$. Thus,  $H_{2n+1}(\Z)$ embeds as a discrete elementary subgroup $\Ga  < PU(n+1,1)$ such that $M_\Ga$ is diffeomorphic to $N\times (0,\infty)$. The partial compactification $\ol{M}_\Ga$ is diffeomorphic to $N\times [0,\infty)$. 
\end{example}

The rest of our examples are nonelementary.

\begin{example}
{\bf Schottky groups.} These are convex-cocompact subgroups $\Ga< G$ isomorphic to free nonabelian groups $F_k$ of finite rank $k$. The 
limit set $\La_\Ga$ is homeomorphic to the Cantor set. Its Hausdorff dimension is positive but can be arbitrarily close to $0$.  Schottky groups are always Stein. 
Every nonelementary discrete subgroup contains a Schottky subgroup. Schottky subgroups can be found via the following procedure. Let $\ga_1,...,\ga_k$ be hyperbolic isometries with pairwise disjoint fixed-point sets. Then there exists $t_0$ such that for each integer $t\ge t_0$, the subgroup generated by $s_1=\ga_1^t,..., s_k=\ga_k^t$ is a Schottky group with the free generating set $s_1,...,s_k$. 
\end{example}

\begin{figure}[tbh]
\centering
\includegraphics[width=100mm]{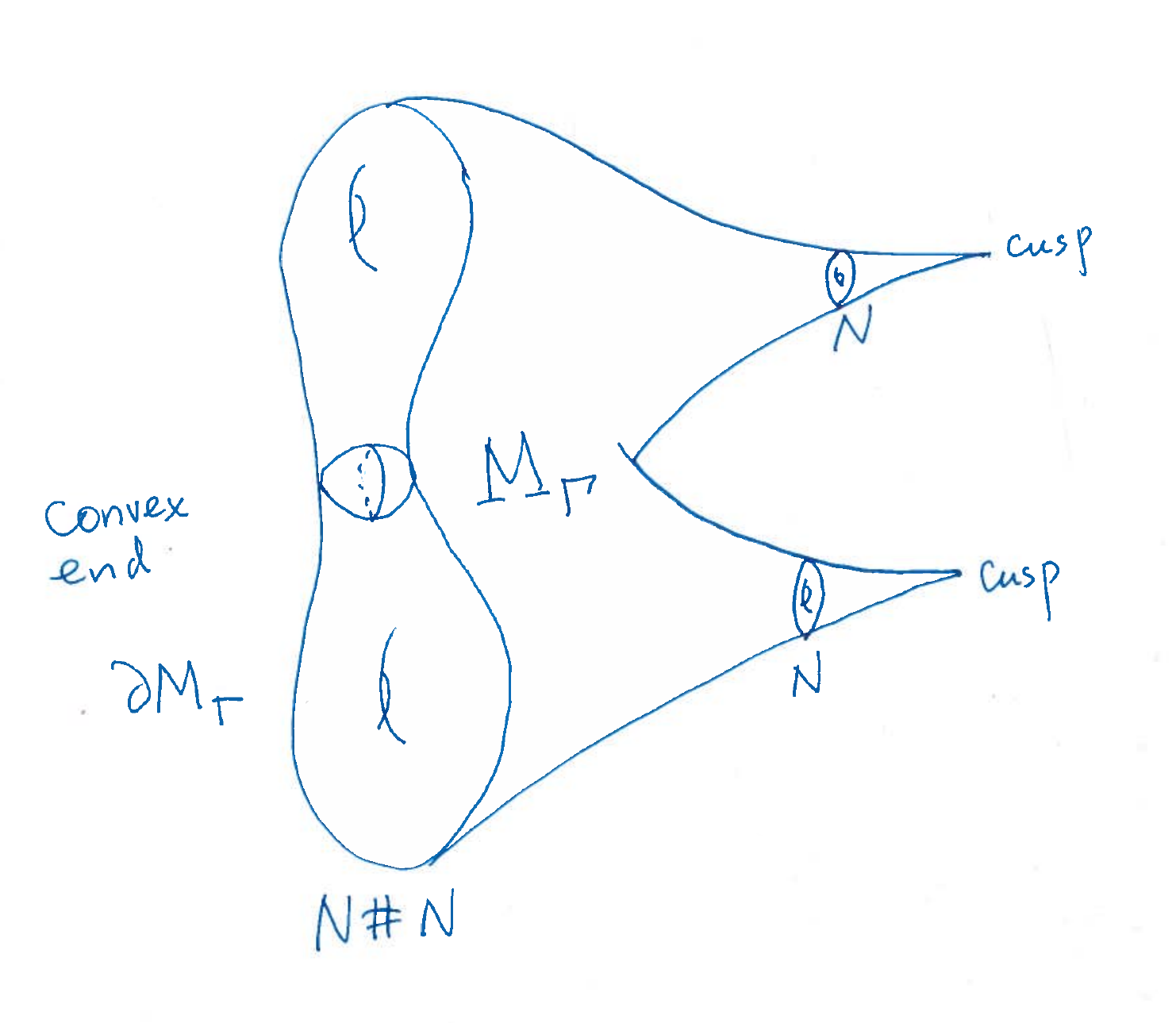}
\caption{Quotient manifold of a Schottky-type group with $k=2$.}
\label{schottky-type.fig}
\end{figure}

\begin{example}\label{ex:Sch-type} 
{\bf Schottky-type groups.} These are geometrically finite subgroups $\Ga< G$ isomorphic to free products of elementary subgroups of $G$, such that 
the limit set $\La_\Ga$ is homeomorphic to the Cantor set.   Schottky-type subgroups can be found via the following procedure. Let $\Ga_1,...,\Ga_k$ be elementary subgroups 
with pairwise disjoint limit sets. Then there exist torsion-free finite-index subgroups $\Ga_i^\ell< \Ga_i, i=1,...,k$, such that the subgroup generated by 
$$
\Ga_1^\ell,...,\Ga_k^\ell 
$$
is Schottky-type and the homomorphism 
$$
\Ga_1^\ell\star...\star\Ga_k^\ell\to \Ga=\<\Ga_1^\ell,...,\Ga_k^\ell\>
$$
sending $\Ga_i^\ell\to \Ga_i^\ell, i=1,...,k$, is an isomorphism.  For instance, suppose that  $\Ga_1,...,\Ga_k$ are integer 
Heisenberg subgroups of $G$. Then $M_\Ga$ has $k$ cuspidal ends 
(diffeomorphic to $N\times (0,\infty)$) and one convex end, with $\partial M_\Ga$ 
diffeomorphic to the $k$-fold connected sum of $N$ with itself, where $N= H_{2n-1}/H_{2n-1}(\Z)$. See Figure \ref{schottky-type.fig}.  
\end{example}

Real and complex Fuchsian groups defined below were introduced by Burns and Shnider in \cite{BS}.

\begin{example}
{\bf Real-Fuchsian subgroups.} Let $\H^2_\R\subset \H^n_\C$ be a real hyperbolic plane in $\H^n_\C$. Let $\Ga< PU(n,1)$ be a geometrically finite subgroup whose limit set is $\geo \H^2_\R$. Then $\Ga$ 
preserves  $\H^2_\R$ and acts on it with quotient of finite area. The quotient surface-orbifold $\Si$ is the convex core of $M_\Ga$. The limit set of $\Ga$ has Hausdorff dimension $1$. 
Assume now that $n=2$, $\Ga$ is torsion-free and $\Sigma$ is compact. Then 
$M_\Ga$ is diffeomorphic to the tangent bundle of $\Sigma$  and is Stein. 
\end{example}

\begin{example}
{\bf Real quasi-Fuchsian subgroups.} Let $\Ga_t, t\in [0,1]$, be a continuous family of discrete convex-cocompact subgroups of $PU(n,1)$ such that $\Ga_0$ is real-Fuchsian but the rest of subgroups $\Ga_t, t>0$, are not.\footnote{Such deformation exist as long as $\Ga_t$ 
is, say, torsion-free. More generally, such deformations exist if $\Ga$ has trivial center and is not isomorphic to a von Dyck group. See e.g. \cite{Weil}.} 
The subgroups $\Ga_t$ ($t>0$) 
are real quasi-Fuchsian subgroups. Their limit sets are topological circles of Hausdorff dimension $>1$. 

Assume that $n=2$, $\Ga$ is torsion-free and $\Sigma$ is compact. Then 
$M_\Ga$ is diffeomorphic to the tangent bundle of $\Sigma$ and is Stein. 
\end{example}

\begin{example}
{\bf Complex Fuchsian subgroups.} Let $\H^1_\C\subset \H^n_\C$ be a complex hyperbolic line in $\H^n_\C$. Let $\Ga< PU(n,1)$ be a geometrically finite subgroup whose limit set is $\geo \H^1_\C$. Then $\Ga$ preserves  $\H^1_\C$ and acts on it with quotient of finite area. The quotient surface-orbifold $\Si$ is the convex core of $M_\Ga$. The limit set of $\Ga$ has Hausdorff dimension $2$. Let $W\subset V=\C^{n+1}$ be the 2-dimensional complex linear subspace such that the projection of $W\cap V_-$ to $\B^n$ equals $\H^1_\C$. The $W^\perp\subset V$ (the complex orthogonal complement  with respect to the form $q$ on $V$) has the property that the form $q$ restricted to $W^\perp$ is positive-definite. The projection $[W^\perp]$ of $W^\perp$ to $\P^n$ is $\Ga$-invariant. The pair $([W], [W^\perp])$ defines a linear holomorphic fibration of $\P^n - [W^\perp]$ over $[W]$: The fiber through $x\in \P^n - [W^\perp]$ is the unique projective hyperplane 
passing through $x$ and intersecting transversally both  $[W]$ and $[W^\perp]$. Restricting to $\B^n$ we obtain a $\Ga$-invariant holomorphic fibration $\B^n\to \H^1_\C$. Projecting to $M_\Ga$ yields a holomorphic orbi-fibration $M_\Ga\to \Si$, whose fibers are biholomorphic to quotients of $\B^{n-1}$ by finite subgroups of $Aut(\B^{n-1})$. 
Assume now that $n=2$, $\Ga$ is torsion-free and $\Sigma$ is compact. Then 
$M_\Ga$ is diffeomorphic to the square root of the tangent bundle of $\Sigma$ (the spin-bundle) and is not Stein (it contains the compact complex curve $\Si$). \end{example}

Convex-cocompact complex Fuchsian groups are locally rigid in the sense that any small deformation of such a group is again complex Fuchsian, \cite{Toledo}.  
The complex Fuchsian examples generalize to the class of geometrically finite subgroups of $PU(n,1)$ whose limit sets are ideal boundaries of $k$-dimensional complex hyperbolic subspaces 
$\H^k_\C\subset \H^n_\C$. The rigidity theorem holds in this case as well, see \cite{GoMi, Corlette0, BI}.

\begin{example}
{\bf Hybrid groups.} One can combine, say, torsion-free, real and complex Fuchsian groups in a variety of ways. For instance, one can form free products of such groups. The nature of the quotient manifolds will depend on the precise way in which the free factors are interacting with each other. For instance, when $n=2$, the boundary of $M_\Ga$ can be either a connected sum, or the {\em toral sum} of certain circle bundles over surfaces.  One can also break real and complex Fuchsian groups into smaller pieces and consider amalgams over $\Z$ of these pieces. As the result, one can get for instance,  circle bundles over surfaces other than the unit tangent bundle and its square root, see \cite{GKL} and \cite{AGG} for more detail. 
\end{example}

\begin{example}
{\bf AGG groups: Anan'in--Grossi--Gusevskii}, \cite{AGG}. These interesting examples of convex-cocompact subgroups of $PU(2,1)$ are isomorphic images of {\em von Dyck groups} $D(2,n,n)$, for $n\in \{10\}\cup [13, 1001]$. None of these subgroups is  complex Fuchsian or real quasi-Fuchsian. According to Proposition \ref{prop:rigid}, 
these subgroups are locally rigid in $PU(2,1)$: Every small deformation is conjugate in $PU(2,1)$ to the original subgroup. 
The limit set is a topological circle but is neither a complex nor a real circle. Fix a (unique up to conjugation) discrete, faithful and isometric  action of  $D(2,n,n)$ on $\H^1_\C$. 
For each embedding $\rho: D(2,n,n)\to \Ga< PU(2,1)$ constructed in section 3.3 of \cite{AGG}, the complex hyperbolic orbifold $M_\Ga$ is diffeomorphic to the total space of an orbifold bundle over the complex 1-dimensional orbifold $\BB=\H^1_\C/D(2,n,n)$ with fibers given by projections to $M_\Ga$ of some complex geodesics in $\H^2_\C$. It follows from the local rigidity of  $\rho$, combined with \cite[Lemma 4.5]{Simpson}, that there exists a $\rho$-equivariant holomorphic map  
$$
\tilde f: \H^1_\C\to \H^2_\C. 
$$
(I owe this observation to Ludmil Katzarkov.)\footnote{I refer the reader to the book \cite{CMP} for a gentle introduction to Simpson's results, discussion of variations of Hodge structures and period domains.} 
Since the orbi-bundle $M=M_\Ga\to \BB$ has holomorphic fibers, it follows that $\tilde f$ descends to a holomorphic map $f: \BB\to M$ which has only positive, zero-dimensional intersections 
with the fibers. Composing with the projection $M\to \BB$, we obtain a self-map $h: \BB\to \BB$ which is a branched covering. Since $\BB$ is a hyperbolic orbifold, it follows that $h=\id$. In other words, $M\to \BB$ admits a holomorphic section. In particular, 
$M$ (and any of its finite manifold-covering spaces, given by Selberg's Lemma) is non-Stein. 
 \end{example}

\begin{example}
{\bf Polygon-groups, J. Granier,} \cite{Gran}. \label{ex:polygon}
The polygon-group $\Ga_{6,3}$ (see Example \ref{ex:poly}) embeds as a convex-cocompact subgroup in $PU(2,1)$ via the reflection representation $\rho_{6,3}$. Thus, the limit set of   $\Ga_{6,3}< PU(2,1)$ is  homeomorphic to the Menger curve. 
\end{example}

Conjecturally, the same holds for all polygon-groups $\Ga_{n,3}, n\ge 6$, cf. \cite{Bourdon, Kapovich2005, DH} for a discussion of isometric actions on real hyperbolic spaces.  

\newpage 

\begin{example}\label{ex:fibrations} 
{\bf Complex-hyperbolic manifolds which are singular fibrations with compact fibers.}  
\end{example} 

\begin{defn}
A  {\em singular Kodaira fibration} is a surjective holomorphic map with connected fibers, $f: M\to B$, between connected complex manifolds/orbifolds, where $0< \dim B< \dim M$. (Usually, it is required that no two generic fibers are biholomorphic to each other, but, in order to simplify the discussion, I will omit this condition.) 
%In this case, we will say that $M$ is {\em fibered}.  
%A complex manifold $M$ is said to be the total space of a 
\end{defn}

Singular Kodaira fibrations need not be locally trivial in holomorphic or even topological sense; a {\em Kodaira fibration} is a holomorphic map  $f: M\to B$ which is a smooth fiber bundle.

In the context of complex hyperbolic manifolds, the first example of a singular Kodaira fibration appeared in Ron Livne's PhD thesis,  \cite{Livne}. Many more examples are now known.  Below we discuss one example which (to my knowledge) first appeared in the work of Hirzebruch, \cite{Hirz}.

\begin{figure}[tbh]
\centering
\includegraphics[width=100mm]{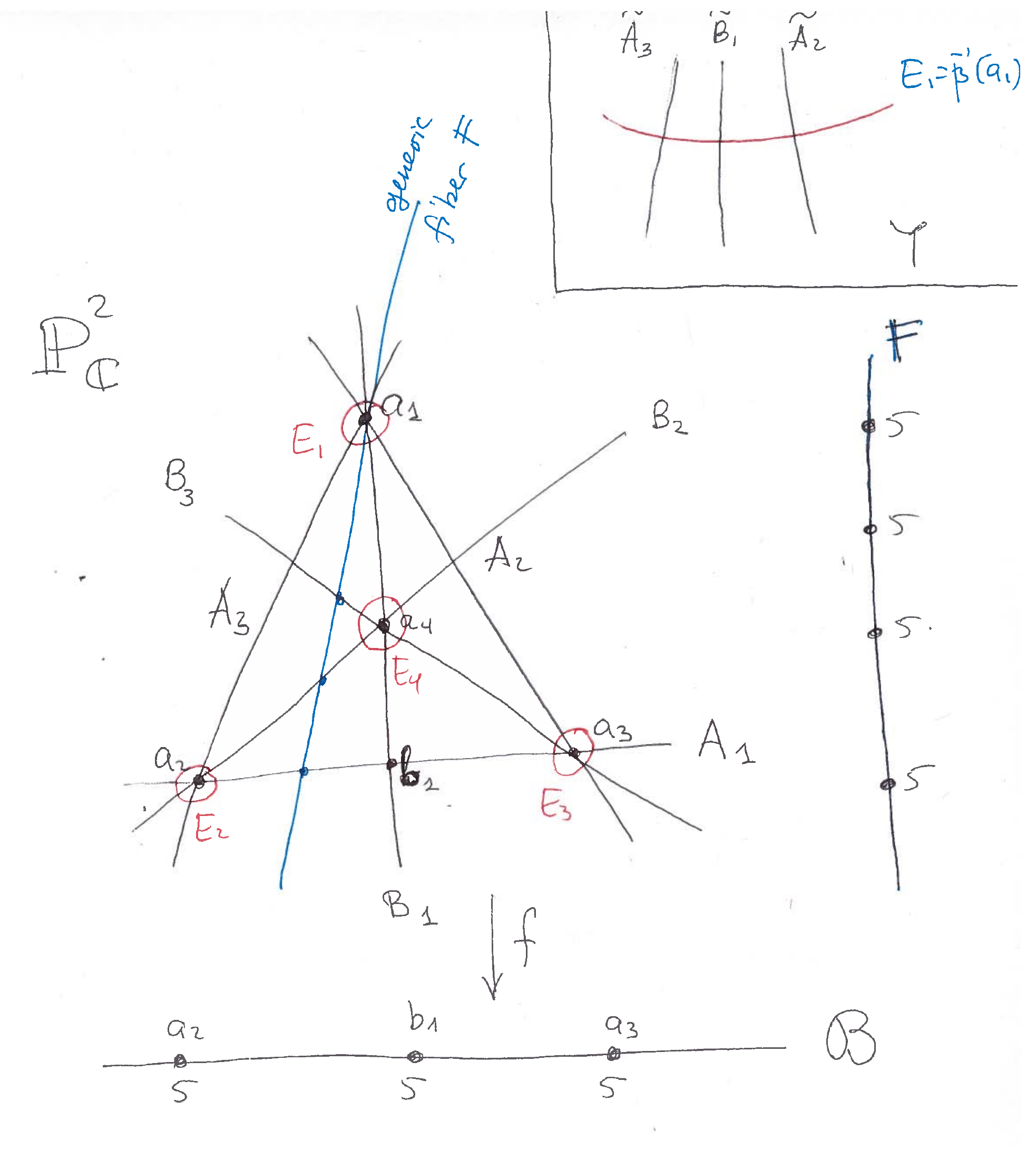}
\caption{Singular Kodaira orbi-fibration}
\label{fib.fig}
\end{figure}

Consider the {\em standard quadrangle} in $\P^2_\C$, which is a configuration $A$ of six lines $A_1, A_2, A_3, B_1, B_2, B_3$ with four triple intersection points $a_1, a_2, a_3, a_4$ and three double intersection points $b_1, b_2, b_3$, see Figure \ref{fib.fig}. Let $Y$ denote the complex surface obtained via blow-up of the four triple intersection points of $A$; let $\beta: Y\to \P^2_\C$ denote the blow-down map. 
Then $Y$ contains a configuration $\tilde{A}$  of eight distinguished smooth rational curves $C_1,...,C_{10}$: 
The four exceptional divisors $E_1,...,E_4$ of the blow-up and six lifts $\tilde{A}_i, \tilde{B}_i, i=1,2,3$, of the original projective lines in the arrangement $A$. The configuration $\tilde{A}$ is a divisor $D$ with simple normal crossings: Any two curves intersect in at most one point and at every intersection point only two curves intersect. Our next goal is to define a {\em complex orbifold} $\O$ with the underlying space $Y$ and the singular/orbifold locus $\Sigma_{\O}$ equals the union of curves in  $\tilde{A}$ (the preimage under $\beta$ of the union of lines in $A$). The local complex orbifold-charts of   $\O$ are defined as follows. 

1. At every point $z\in \O - \Sigma_{\O}$ the local chart is given by the restriction of $\beta$ to a suitable neighborhood of $z$. 

2. At every point $z\in  \Sigma_{\O}$ which is not a (double) intersection point of the divisor $D$ but $z\in C_i$, $i=1,...,10$, the local chart is the holomorphic $5$-fold branched covering over a suitable neighborhood of $z$, ramified over $C_i$. 

3. Suppose that $z$ is an intersection point of $D$, $z\in C_i\cap C_j$, $i\ne j$. Choose local holomorphic coordinates at $z$ where $C_i, C_j$ correspond to the coordinate lines in $\C^2$ and $z$ corresponds to the origin; $\C^2=\C\times \C$. Each factor $\C$ in this product decomposition is biholomorphic to the quotient $\C/\Z_5$, with $\Z_5$ acting linearly  on $\C$. Thus, a small neighborhood $U$ of $z$ in $Y$ is biholomorphic to the quotient of the bi-disk, $\Delta^2/\Z_5^2$. This yields the local orbifold-chart at $z$, $\Delta^2\to \Delta^2/\Z_5^2\cong U$. 

The result is a complex orbifold $\O$ with the underlying space $Y$. Hirzebruch then proves that the orbifold $\O$ is biholomorphic to the orbifold-quotient $M_\Ga=\B^2/\Ga$ of the complex 2-ball, by appealing to Yau's Uniformization Theorem, \cite{Yau}: He verifies that the orbifold $\O$ admits a finite holomorphic orbifold-covering $M\to \O$, where $M$ is a complex surface of general type satisfying the equality of characteristic classes $3c_2=c^2_1$; equivalently, $3\sigma(M)=\chi(M)$, where $\sigma$ is the signature and $\chi$ is the Euler characteristic.  Yau's theorem implies that $M$ admits a K\"ahler metric of constant bisectional curvature $-1$, i.e. is a ball-quotient. Mostow Rigidity Theorem then implies that $\O$ is a complex hyperbolic orbifold as well. A bit more streamlined version of this argument was later developed by Barthel--Hirzebruch--Hofer, \cite{BHH}, and  Holzapfel, \cite{Holz}, who defined orbifold-characteristic classes directly computable from  lines arrangement $A$ in $\P^2_\C$ (as well as $\P^1_\C \times \P^1_\C$) and the orbifold-ramification numbers assigned to rational curves in the corresponding post-blow-up divisor. 

I next describe a singular Kodaira orbi-fibration on $\O$. Pick one of the triple intersection points, say, $a_1$, 
of the arrangement $A$ and let $A_1$ be a line in $A$ not passing through $a_1$. Consider the pencil of projective lines passing through $a_1$. This pencil defines a (nonsingular) holomorphic fibration of $\P^2_\C -\{a_1\}$ with the base $A_1$; the fibration map sends $z\in \P^2_\C -\{a_1\}$
to the point of intersection of the line $z a_1$ with the line $A_1$. This fibration becomes a holomorphic map $f: Y\to \tilde{A}_1$ when we lift it to $Y$. Some fibers of $f$ are, however, singular: These are the three singular fibers corresponding to the lifts of the  three lines $A_2, A_3, B_1$ passing through $a_1$ and other points of triple intersection of $A$: $a_2, a_3, a_4$. The corresponding fibers are reducible rational curves (with the extra components corresponding to the exceptional divisors $E_2, E_3, E_4$). The line $A_1$ has an orbifold structure induced from $\O$: The corresponding  orbifold ${\mathcal B}$ has three singular points $a_2, a_3, b_1$, with the local isotropy groups $\Z_5$ for each of them. 
The map $f$ defined above respects the orbifold structure of $\O$ and ${\mathcal B}$ and, hence, we obtain a singular Kodaira orbi-fibration $f: \O \to {\mathcal B}$. This fibration is nonsingular away from the preimages of the points $a_2, b_1, a_3$, with the generic fiber(s) $\F$ diffeomorphic to the orbifold with the underlying space $\P^1_{\C}$ and four singular points of the order $5$.

The restriction of $f$ to $\O'=f^{-1}(\{a_2, b_1, a_3\})$ is a nonsingular Kodaira fibration, i.e. a smooth (orbifold) fiber bundle; accordingly, $\pi_1(\F)$ embeds as a normal subgroup 
 in $\pi_1(\O')$. Since the inclusion $\O'\to \O$ induces an epimorphism of fundamental groups $\pi_1(\O')\to \pi_1(O)=\Ga$,  the image $N$ of   $\pi_1(\F)$ in 
$\pi_1(O)=\Ga$ is a normal finitely-generated subgroup $N\triangleleft \Ga$. By passing to the universal covering of $\BB$, we obtain a  holomorphic map $h: \H^2_\C/N \to \H^1_\C$. The fibers of this map are compact and, generically, diffeomorphic to $\F$. The map $h$ has infinitely many critical values in  $\H^1_\C$ which break into finitely many $\pi_1(\BB)$-orbits and accumulate to the entire circle $\geo \H^1_\C$.  Lifting $h$ further to an $N$-invariant holomorphic function $\H^2_\C\to \H^1_\C$ and extending this function to a 
 measurable $N$-invariant nonconstant function $S^3=\geo \H^2_\C\to S^1$, we conclude that the action of $N$ on $S^3$ is non-ergodic. 

\medskip
The group $\Ga$ in the above example is a special case of: 

\begin{example}
{\bf Arithmetic lattices of  simplest type.} Let $K$ be a totally real number field, i.e. a finite extension of $\Q$ such that the image of every embedding $K\to \C$ lies in $\R$. 
Take an imaginary quadratic extension $L/K$, i.e. an extension which does not embed in $\R$. Since $K$ is totally-real and $L$ is its imaginary extension, all embedding 
$L\to \C$ come in complex conjugate pairs:
$$
\si_1, \bar\si_1,..., \si_k, \bar\si_k. 
$$
Next, take a hermitian quadratic form in $n+1$ variables 
$$
\varphi(z, \bar{z})=\sum_{p,q=1}^{n+1} a_{pq} z_p \bar{z}_q
$$
with coefficients in $L$. We require the forms $\varphi^{\si_1}, \varphi^{\si_2}$ to have the signature $(n,1)$  
and the forms $\varphi^{\si_j}, \varphi^{\bar\si_j}$ to be definite for the rest of the embeddings.  
We will identify $L$ with $\si_1(L)$, so $\si_1=\id$. Let $SU(\varphi)$ denote the group of special unitary automorphisms of the form $\varphi$ on 
$L^{n+1}$. The embedding $\si_1$ defines a homomorphism $SU(\varphi)\to SU(n,1)$ with relatively compact kernel. 

A subgroup $\Ga$ of $SU(n,1)$ is said to be an arithmetic lattice of the {\em simplest type} (or {\em of type I}) if it is commensurable\footnote{I.e. the intersection of the two groups has finite index in both} to $SU(\varphi, \O_L)= SU(\varphi)\cap SL(n+1, \O_L)$, where $\O_L$ is the ring of integers of $L$. For every such $\Ga$ the quotient $\H^n_\C/\Ga$ has finite volume. I refer to \cite{McReynolds} 
and \cite{Paupert}  for more detail on arithmetic subgroups of $SU(n,1)$. 
\end{example}

It is known that every arithmetic lattice $\Ga$ of the simplest type contains a finite index congruence-subgroup $\Ga'$ with infinite abelianization, \cite{Kazhdan} (see also \cite{Wallach}). Equivalently, the quotient-space $\B^n/\Ga'$ has positive 1st Betti number. 
In contrast, Rogawski,  \cite{Rogawski}, proved that for {\em type I\!I} arithmetic lattices in $SU(2,1)$, every congruence-subgroup has finite abelianization. It is unknown if such a lattice contain finite index subgroups with infinite abelianization. Furthermore, certain classes of non-arithmetic lattices in $SU(2,1)$ (the ones violating the integrality condition for arithmetic groups) are proven to have positive virtual first Betti number by the work of S.-K. Yeung, \cite{Yeung}. 

We now discuss the existence of  (singular) Kodaira fibrations of {\em compact} complex hyperbolic manifolds 
$M=\H^n_\C/\Ga$. % to the existence of , assuming that $M$ is compact and $\Ga$ is torsion-free. 

1. Suppose that $b_1(M)>0$. Since $M$ is K\"ahler, $b_1(M)$ is even, hence, there exists an epimorphism $\phi: \Ga\to \Z^2$. If the kernel of $\phi$ is not finitely-generated, then, according to a theorem of Delzant, \cite{Delzant}, the manifold $M$ admits a singular Kodaira fibration over a 1-dimensional complex hyperbolic orbifold. 

2. If $M$ and $B$ are both complex hyperbolic, then there are no (nonsingular) Kodaira fibrations $M\to B$:  It was first proven in the case when $M$ is a surface by Liu, \cite{Liu}, and then generalized to arbitrary dimensions by Koziarz and Mok, \cite{KM}. They also prove nonexistence of Kodaira fibrations $M\to B$ when $\dim(B)\ge 2$ and $M$ merely has finite volume. 
Furthermore, if $M$ is 2-dimensional, for every singular Kodaira fibration $M\to B$, 
the kernel of the homomorphism $\pi_1(M)=\Ga\to \pi_1(B)$ is finitely generated but is not finitely-presentable, \cite{Kapovich1998, HK}. 

\begin{question}
Is there a discrete  subgroup $\Ga< PU(2,1)$ isomorphic to the fundamental group of a compact real hyperbolic surface, such that $M=M_\Ga$ admits a Kodaira fibration  (with compact fibers) 
$M\to \H^1_\C$? Is there a singular Kodaira fibration (with compact fibers) $\H^2_\C/\Ga \to \H^1_\C$ which has only finitely many singular fibers? 
\end{question}

\section{Complex Kleinian groups and function theory on complex hyperbolic manifolds}\label{sec:interactions} 

In this section we discuss some interesting interactions between the general theory of holomorphic functions on complex manifolds 
(which I review in Section \ref{sec:CA}) 
and geometry/topology  of complex Kleinian groups.

\begin{prop}
Suppose that $\Ga< PU(n,1)$ is a discrete, torsion-free subgroup such that $M=M_\Ga$ admits a surjective holomorphic map 
with compact fibers $f: M\to B$, where $B$ is a complex manifold satisfying $\dim(B)< n$. Then $\Om_\Ga=\emptyset$. 
In particular, $M$ cannot have convex ends.   
\end{prop}
\proof Suppose, to the contrary, that $\Om_\Ga\ne \emptyset$. Then $Core_\eta(M)$ is a proper submanifold (with boundary) in $M$.  Since $\H^n_\C$ is strictly negatively curved, 
the nearest-point projection $\Pi: \H^n_\C\to \hull_\eta(\La_\Ga)$ is strictly contracting away from $\hull_\eta \La_\Ga$. 
By the $\Ga$-equivariance, $\Pi$ descends to a strictly contracting 
projection $\pi: M \to Core_\eta (M)$.   Therefore, if $Y$ is a compact complex subvariety in $M$ of positive dimension not contained 
$Core_\eta (M)$ then $\pi(Y)$ has volume strictly smaller than that of $Y$.  This is a contradiction since  
$\pi: Y\to \pi(Y)$ is homotopic to the identity inclusion map $\id_Y: Y\to M$ and compact complex subvarieties in K\"ahler manifolds are volume-minimizers in their homology classes. 
Taking a generic fiber $Y$ of  $f: M\to B$ through a point $x\in M- \hull_\eta(\La_\Ga)$ concludes the proof. \qed 

\medskip
We next discuss geometry and topology of quotient-orbifolds $M_\Ga$, primarily for convex-cocompact subgroups $\Ga< PU(n,1)$.

A classical example of a complex submanifold with strictly Levi-convex boundary is a closed unit ball $\ol{\B^n}$ in $\C^n$. Suppose that $\Ga< Aut(\B^n)$ is a discrete torsion-free subgroup of the group of holomorphic automorphisms of $\B^n$ with (nonempty) domain of discontinuity $\Omega\subset \partial \B^n$. The quotient 
$$
\ol{M}_\Ga= (\B^n \cup \Om)/\Ga
$$
is a smooth submanifold with strictly Levi-convex boundary boundary in the complex manifold $\check{\Om}_\Ga/\Ga$ (see \eqref{eq:thull}). 
Thus, we conclude:

\begin{lemma}
If $\ol{M}_\Ga=(\B^n \cup \Om)/\Ga$ has compact boundary, then $M$ is strongly pseudoconvex. 
\end{lemma}

Consequently: 

\begin{thm}\label{thm:conn} 
Let $\Ga< PU(n,1)$, $n\ge 2$, be a convex-cocompact discrete subgroup. Then $\partial M_\Ga$ is connected. 
\end{thm}
\proof Since $\Ga$ is convex-cocompact, it is also finitely generated. Hence, by Selberg's Lemma, 
 the orbifold $M_\Ga$ is very good. Therefore, it suffices to consider the case when $\Ga$ is torsion-free, i.e. $M_\Ga$ is a complex $n$-manifold. Since $\ol{M}_\Ga$ 
is strongly pseudoconvex, connectedness of its boundary is an immediate consequence of Theorem \ref{thm:Kohn-Rossi}. \qed

\begin{thm}
Let $\Ga< PU(n,1)$, $n\ge 2$, be a convex-cocompact discrete subgroup which is not a lattice, i.e. $\Om_\Ga\ne \emptyset$. Then $\dim (\La_\Ga) \le 2n-3$, equivalently, $\cd_\Q(\Ga)\le 2n-2$.
\end{thm}
\proof As before, it suffices to consider the case of torsion-free groups $\Ga$. 
According to Corollary \ref{cor:CW}, $M_\Ga$ is homotopy-equivalent to a CW complex of dimension $\le 2n-2$. It follows that  $\cd_\Q(\Ga)\le 2n-2$ and, by the Bestvina-Mess theorem, 
$\dim (\geo \Ga) \le 2n-3$. Since $\geo \Ga$ is homeomorphic to $\La_\Ga$, $\dim (\La_\Ga) \le 2n-3$ as well. \qed 

\medskip
In particular, $\La_\Ga$ does not separate $S^{2n-1}$ (even locally) and, hence, $\Om_\Ga$ is connected, which gives another proof of the fact that $\partial M_\Ga$ is connected. 

Specializing to the case $n=2$, we obtain: If $\Ga< PU(2,1)$ (for simplicity, torsion-free) is convex-cocompact and is not a lattice, then $\La_\Ga$ is at most 1-dimensional. In particular, according to \cite{KK}, $\Ga$ admits 
an iterated amalgam decomposition over trivial and  cyclic subgroups, so that the terminal groups are either 
cyclic, or isomorphic to Fuchsian groups (and the limit set is a topological circle) or groups whose limit sets are 
Sierpinski carpets or Menger curves.

\begin{thm}
Suppose that $\Ga$ is torsion-free convex cocompact, $n>1$ and $M_\Ga$ contains no compact complex subvarieties of positive dimension. 
Then $M_\Ga$ is Stein.  
\end{thm}
\proof This is an immediate consequence of Theorem \ref{thm:rossi}. \qed 

\medskip
One way to prove that $M_\Ga$ contains no compact complex subvarieties of positive dimension is to argue that $\Ga=\pi_1(M)$ is free: This implies that $H_i(M_\Ga)=0, i\ge 2$, but, since $M_\Ga$ is K\"ahler, every compact complex $k$-dimensional subvariety of $M_\Ga$ would define a nonzero $2k$-dimensional homology class. For instance, if $\Ga$ is convex-cocompact, 
$\delta_\Ga<1$ then $\dim \La_\Ga\le \dim_H(\La_\Ga)<1$, which implies that   $\dim \La_\Ga=0$ and, hence, $\Ga$ is a virtually free group. However, even when $H_2(M)\ne 0$, one can still, sometimes, prove that $M_\Ga$ contains no compact complex curves. For instance, let $L\to M_\Ga$ be the canonical line bundle. If $C\subset M_\Ga$ is an (even singular) complex curve, 
the pull-back of $L$ to $C$ has nonzero 1st Chern class. Assuming that $H_2(M_\Ga)\cong \Z$ (e.g. if $\Ga$ is isomorphic to the fundamental group of a compact Riemann surface), if 
the 1st Chern class of $L$ evaluated on the generator of  $H_2(M_\Ga)$ is zero, then $M_\Ga$ contains no complex curves. This argument applies in the case of real Fuchsian groups and their quasi-Fuchsian deformations. 

Observe that if $\Ga< PU(2,1)$ is a complex Fuchsian group, then $\dim_H(\La_\Ga)=2$. 

\begin{thm}
[S. Dey, M. Kapovich, \cite{DeK}]\label{thm:DK} 
 If $\Ga< PU(n,1)$ is discrete, torsion-free and $M_\Ga$ contains a compact complex subvariety of positive dimension, then $\delta_\Ga\ge 2$. 
\end{thm}

\begin{cor}
Suppose that $\Ga< PU(n,1)$ is  torsion-free, convex-cocompact and $\delta_\Ga<2$, then $M_\Ga$ is Stein. 
\end{cor}

%Can we conclude that $\Om$ is connected and its boundary does not locally separate; has topological codimension $\ge 2$? 

\medskip
{\bf Burns' Theorem.} We now drop the convex-cocompactness assumption and consider general  discrete, torsion-free subgroups $\Ga< PU(n,1)$. 
Theorem \ref{thm:conn} has the following striking generalization. It was first stated by Dan Burns, who, as it appears, never published a proof; a published proof is due to Napier and Ramachandran, \cite[Theorem 4.2]{NR}:

\begin{thm}\label{thm:burns} 
Suppose that $n\ge 3$, $\Ga< PU(n,1)$ is discrete, torsion-free and $\partial M_\Ga$ has at least one compact component  $S$. Then:

1. $\D M_\Ga=S$.  

2. $\Ga$ is geometrically finite. 
\end{thm}

%Maybe also that $\Om=\Om_0$...

A good example illustrating this theorem is that of a Schottky-type group (Example \ref{ex:Sch-type}), where the limit set is totally disconnected, the quotient manifold $\Om_\Ga/\Ga$ is compact and $M_\Ga$ has $k$ cusps. In particular, $\ol{M}_\Ga$ is noncompact in this example. 

It is unknown if Burns' theorem holds for $n=2$, but Mohan Ramachandran proved the following:

\begin{thm}\label{thm:mohan}
Suppose that $\Ga< PU(2,1)$ is discrete, torsion-free, the injectivity radius of $M_\Ga$ is bounded away from zero, and $\partial M_\Ga$ has at least one compact component. Then 
$\Ga$ is convex-cocompact. 
\end{thm}

The proof of this theorem is given in Appendix G.

\section{Conjectures and questions}\label{sec:conjectures} 

In this section I collect some conjectures and questions, in addition to those scattered throughout these notes. 

\medskip 
The first conjecture is a generalization of Burns' theorem, Theorem \ref{thm:burns}: 

\begin{conjecture}
Suppose that $\Ga< PU(n,1), n\ge 2$, is such that for $M=M_\Ga$ the thick part  $M_{[\eps,\infty)}$ has a convex end. Then $\Ga$ is geometrically finite and $\Om_\Ga$ is connected.  
\end{conjecture}

The next two conjectures are motivated by Theorem  \ref{thm:DK}: 

\begin{conjecture}
If $\Ga< PU(n,1)$ is discrete, torsion-free, non-Stein subgroup, then $\Ga$ is a complex Fuchsian group. 
\end{conjecture}

\begin{conjecture}
If $\Ga< PU(n,1)$ is discrete, torsion-free and $\delta_\Ga<2k$, then $M_\Ga$ cannot contain a compact complex subvariety of dimension $k$. 
\end{conjecture}

%\begin{conjecture}
%Suppose that $\Ga< Aut(\B^n)$ is a convex-cocompact torsion-free subgroup and $\delta_\Ga<2$. Then $\B^n/\Ga$ is Stein.   
%In the case of equality, $\Ga$ is a complex Fuchsian subgroup: a subgroup acting cocompactly on a complex 1-disk in $\B^n$.
%\end{conjecture}

\begin{conjecture}
[Chengbo Yue's Gap Conjecture, \cite{Yue}] Suppose that $\Ga< G= Aut(\B^n)$ is a convex-cocompact torsion-free subgroup.  Then either $\Ga$ is a uniform lattice in $G$ (and, thus, $\delta_\Ga=2n$) 
or $\delta_\Ga\le 2n-1$. 
\end{conjecture}

Note that the two other conjectures about nonelementary convex-cocompact subgroups $\Ga< PU(n,1)$ 
made in the introduction to \cite{Yue} fail already in dimension $n=2$:

(a) The inequality $\dim_H \La_\Ga > n-1$ does not imply that $M_\Ga$ is non-Stein. For instance, a real-hyperbolic quasifuchsian subgroup of $PU(2,1)$ serves as an example. 

(b) Even if $M_\Ga$ is non-Stein, a compact complex curve  in $M_\Ga$ need not be a finite union of totally geodesic complex curves, as it is shown by the AGG-examples.

%\begin{conjecture}
%Suppose that $\Ga< PU(2,1)$ is such that $M_\Ga$ has a convex end. Then $\Ga$ is geometrically finite and $\Om_\Ga$ is connected.  
%\end{conjecture}

\begin{problem}
1. Investigate which polygon-groups embed discretely in $PU(2,1)$. 

2. Is there a convex-cocompact subgroup of $PU(2,1)$ with the limit set  homeomorphic to the Sierpinski carpet?  
\end{problem}

While ``most'' compact 3-dimensional manifolds are hyperbolic, very few examples of hyperbolic 3-manifolds which are of the form $\Om_\Ga/\Ga$, $\Ga< PU(2,1)$ are known, see 
the book by Richard Schwartz \cite{Sch} for further discussion.

\begin{conjecture}
The Menger curve limit set in Example \ref{ex:polygon} is ``unknotted'' in $S^3$, i.e. is ambient-isotopic to the standard Menger curve 
${\mathcal M}\subset \R^3\subset S^3= \R^3\cup \{\infty\}$. 
Furthermore, in this example, 
the quotient 3-dimensional manifold $\Om_\Ga/\Ga$ is hyperbolic.\footnote{It suffices to show that $\Om_\Ga/\Ga$ contains no incompressible tori, which is closely related to the unknottedness 
problem of the Menger-curve limit set.} 
\end{conjecture}

\begin{problem}
Prove the existence of discrete geometrically infinite subgroups of $PU(2,1)$ isomorphic to the fundamental groups of compact surfaces.\footnote{Cf. section 11.4 in \cite{Kapovich2008}.}  
\end{problem}

Note that such subgroups do not exist in $PU(1,1)$ but abound in $PO(3,1)$. Furthermore, the only known examples of finitely generated geometrically infinite subgroups of $PU(2,1)$ 
come from singular Kodaira fibrations and are not finitely-presentable, see Example \ref{ex:fibrations}.

\bigskip 
The conjectures and questions appearing above, deal with discrete subgroups $\Ga$ of $PU(n,1)$ which are not lattices, i.e. the $\H^n_\C/\Ga$ has infinite volume. 
Below, I discuss two problems regarding  lattices. 

\medskip
{\bf Arithmeticity.} The most famous open problem regarding lattices in $PU(n,1)$ deals with the existence problem of nonarithmetic subgroups and was first raised in Margulis' ICM address \cite{Margulis}. It is known (due to the work of Margulis \cite{Margulis-book}, Corlette \cite{Corlette-NA}, Gromov--Schoen \cite{GS}, and Gromov--Piatetski-Shapiro \cite{GPS}) that:

(a) For each $n$, the Lie group $SO(n,1)$ contains non-arithmetic lattices.

(b) For every simple noncompact connected linear Lie group $G$ which is not locally isomorphic to $SO(n,1)$ and $SU(n,1)$, every lattice $\Ga< G$ is arithmetic. 

This leaves out the series of Lie groups $PU(n,1)$, $n\ge 2$. Currently, primarily due to the work of Deligne and Mostow, see \cite{DM}, there are known examples of nonarithmetic lattices in 
$PU(2,1)$ and $PU(3,1)$. Loosely speaking there are three approaches to constructing nonarithmetic lattices:

(a) As monodromy groups of some linear holomorphic ODEs,  see  \cite{DM, CHL}, as well as \cite{Th} for 
a geometric interpretation. 

(b) By constructing the corresponding complex hyperbolic orbifolds $M_\Ga$ whose underlying space is a blown-up $\P^n$, see \cite{BHH, Shv, CHL, Deraux}

(c) By constructing a Dirichlet fundamental domain of $\Ga$ in $\H^2_\C$, see \cite{Mostow1980, DPP}.    

\medskip 
But using these techniques becomes increasingly difficult (or even impossible) as the dimension $n$ increases, which means that different approaches are needed. 

\begin{conjecture}
For each $n$, $PU(n,1)$ contains a nonarithmetic lattice. 
\end{conjecture}

By analogy with the construction of non-arithmetic lattices in \cite{GPS}, one can hope for a similar ``hybrid'' construction of nonarithmtic lattices in $PU(n,1)$, leading to a conjecture due to Bruce Hunt:

\begin{conjecture}
For every $n\ge 2$, there exists a quadruple of arithmetic lattices $\Ga_1, \Ga_2< SU(n-1,1)$ and $\Ga_3< SU(n-2,1)$ such that:

(1)  $\Ga_3$ is isomorphic to subgroups in  $\Ga_1, \Ga_2$; hence, we obtain an amalgam $\Ga_0= \Ga_1\star_{\Ga_3} \Ga_2$. 

(2) There exists an epimorphism  $\rho: \Ga_0\to \Ga< SU(n,1)$ injective on $\Ga_1, \Ga_2$, whose image is a nonarithmetic lattice 
$\Ga< SU(n,1)$. 
\end{conjecture}

Unlike \cite{GPS}, where  nonarithmetic lattices in $SO(n,1)$ were constructed via a similar process, with an {\em isomorphism} $\Ga_0\to \Ga< SO(n,1)$, 
in the complex hyperbolic setting there is no hope for an injective homomorphism $\rho$ (a lattice in $SU(n,1)$ cannot be isomorphic to an amalgam  $\Ga_0$ as above). 

\medskip
{\bf Nonexistence of reflection lattices.} The known examples of nonarithmetic lattices $\Ga$ in $PU(n,1), n=2, 3$, are all commensurable to {\em complex reflection subgroups}, i.e. discrete 
subgroups of $PU(n,1)$ generated by complex reflections. Furthermore, up to commensuration, the underlying spaces of their quotient orbifolds $M_\Ga=\H^n_\C/\Ga$ 
are {\em rational projective varieties}. 

\begin{conjecture}
There exists $N$ such that for all $n\ge N$ the following holds:

1. If $\Ga< PU(n,1)$ is a lattice then $\Ga$ cannot be a reflection subgroup. 

2. If $\Ga< PU(n,1)$ is a lattice then the underlying space of the orbifold $M_\Ga$ cannot be a rational algebraic variety. More ambitiously, it has to be a variety of  general type. 
\end{conjecture}

The motivation for this conjecture comes from theorems due to Vinberg, \cite{Vinberg}, and Prokhorov, \cite{Prokhorov}, establishing nonexistence of reflection 
lattices in $PO(n,1)$, when $n$ is sufficiently large.

\section{Appendix A. Horofunction compactification} 
\label{sec:horoboundary}

A metric space $(Y,d)$ is called {\em geodesic} if any two points $x, y$ in $X$ are connected by a geodesic segment, denotes $xy$. (This notation is a bit ambiguous since in many cases such a segment is non-unique.)  A {\em geodesic triangle}, denoted $xyz$, in a metric space $(X,d)$ is a set of three geodesic segments $xy, yz, zx$ connecting cyclically the points $x, y, z$, the {\em vertices} of the triangle; the segments $xy, yz, zx$ are the {\em edges} of the triangle. Thus, geodesic triangles are 1-dimensional objects.

Let $(Y,d)$ be a locally compact geodesic metric space. For each $y\in Y$ define the 1-Lipschitz function 
$d_y= d(y, \cdot)$ on $Y$.  This leads to the {\em Kuratowski embedding} $\kappa: Y\to C(Y)=C(Y,\R)$, $y\mapsto d_y$.  We let $\R\subset C(Y)$ denote the linear subspace of constant functions.  Composing the embedding $\kappa$ with the projection $C(Y)\to C(Y)/\R$ (where $\R$ acts additively on $C(Y)$) 
we obtain the {\em Kuratowski  embedding} of $Y$,
$$
Y\embed  C(Y)/\R. 
$$ 
Then $\ol{Y}$, the closure of $Y$ in ${C(Y)}/\R$, is the {\em horofunction compactification} of $Y$. 
\index{horofunction compactification}
Functions representing points in $\geo Y= \ol{Y}- Y$ are the {\em horofunctions} on $Y$. In other words, horofunctions on $Y$ are limits (uniform on compacts in $Y$) of sequences of normalized distance functions $d_{y_i} - d_{y_i}(o)$, where $y_i\in Y$ are divergent sequences  in $Y$. Each geodesic ray $r(t)$ in $Y$ determines a horofunction in $Y$ called a {\em Busemann function} $b_r$, which is the subsequential limit
$$\index{Busemann function}
\lim_{i\to\infty} d_{r(i)} - d_{r(i)}(o). 
$$
If $Y$ is  a Hadamard manifold, then each limit as above exists (without passing to a subsequence). Furthermore, each horofunction is a Busemann function.  This yields a topological identification of the visual compactification of $Y$ and its horofunction compactification. Level sets of Busemann functions are called {\em horospheres} \index{horosphere}
in $X$. The point $r(\infty)\in \geo Y$ is the {\em center} of the horosphere $\{b_r=c\}$.  Sublevel sets 
$\{ b_r< c\}$ are called {\em horoballs}. The point $r(\infty)$ represented by the ray $r$ is the {\em center} of the corresponding horospheres/horoballs.

\section{Appendix B: Two classical Peano continua} 

A {\em Peano continuum} is a compact connected 
and locally path-connected metrizable topological space. 
We will need two examples of 1-dimensional Peano continua.  
Both are obtained via a procedure similar to the construction of the ``ternary'' Cantor set.

\medskip
{\bf Sierpinski carpet.} Let $I=[0,1]$ denote the unit interval. Start with the unit square $Q_0=I^2\subset \R^2$. Divide $I$ in three congruent subintervals and, accordingly, divide $I^2$ in 9 congruent subsquares. Remove the interior of the ``middle subsquare'', the one disjoint from the boundary of $Q$. Call the result $Q_1$. Now, repeat this procedure for each of the remaining 8 subsquares in $Q_1$, to obtain a planar region $Q_2$, etc. The {\em standard Sierpinski carpet} in $\R^2$ is the intersection
$$
{\mathcal S}:= \bigcap_{i=0}^\infty Q_i. 
$$ 

\medskip
{\bf Menger curve.} Consider the unit cube $C=I^3\subset \R^3$. Let $\pi_i, i=1, 2, 3$ denote the orthogonal projections of $\R^3$ to the coordinate hyperplanes $P_i, i=1,2,3$, in $\R^3$. In all three planes we take the Sierpinski carpets ${\mathcal S}_i\subset P_i$, 
constructed from the unit squares $Q_i=C\cap P_i$, $i=1,2,3$. Then the {\em standard Menger curve} in $\R^3$ is defined as
$$
{\mathcal M}:= \bigcap_{i=1}^3 \pi_i^{-1}( {\mathcal S}_i). 
$$

\section{Appendix C: Gromov-hyperbolic spaces and groups}

A geodesic metric space  $(X,d)$ is called {\em $\delta$-hyperbolic} if every geodesic triangle $xyz$ in $X$ is $\delta$-slim, i.e. every edge of $xyz$ is contained in the closed $\delta$-neighborhood of the union of the other two edges.  A  geodesic metric space  is called {\em Gromov-hyperbolic} if  it is $\delta$-hyperbolic for some $\delta<\infty$. 

Examples of Gromov-hyperbolic spaces are {\em strictly negatively curved Hadamard manifolds}: If $X$ is a Hadamard manifold of sectional curvature $\le -1$ then $X$ is $\delta_0$-hyperbolic with $\delta_0=\arccosh(\sqrt{2})$. 

Let $\Ga$ be a group with finite generating set $S$. Given $S$, one defines the {\em Cayley graph} $C_{\Ga,S}$. This graph is connected and $\Ga$ acts on it with finite quotient (the quotient graph has a single vertex and $\card(S)$ edges).  The graph $C_{\Ga,S}$ has a graph-metric, where every edge has unit length. 

\begin{defn}
A finitely generated group $\Ga$ is called {\em Gromov-hyperbolic} or simply {\em hyperbolic} if one (equivalently, every) Cayley graph of $\Ga$ is a Gromov-hyperbolic metric space.  
\end{defn}

The {\em Gromov boundary} $\geo \Ga$ of $\Ga$ is the horoboundary of one (any) Cayley graph of $\Ga$: Gromov boundaries corresponding to different Cayley graphs are equivariantly homeomorphic. 

\medskip 
Examples of hyperbolic groups are given by:

\begin{example}
Let $X$ be a strictly negatively curved Hadamard manifold, $Y\subset X$ is a closed convex subset and $\Ga< \Isom(X)$ acts properly discontinuously and cocompactly on  $Y$. Then $\Ga$ is hyperbolic and the ideal boundary $\geo Y$ of $Y$ is equivariantly homeomorphic to the Gromov boundary of $\Ga$. 
\end{example}

In particular, every convex-cocompact discrete subgroup $\Ga< \Isom(X)$ is Gromov-hyperbolic and $\geo \Ga$ is equivariantly homeomorphic to the limit set of $\Ga$. 

\medskip 
Cohomological dimension (with respect to the Chech cohomology) of the Gromov boundary of a hyperbolic group is closely related to the rational cohomological dimension of $\Ga$ itself: 

\begin{thm}
[Bestvina--Mess, \cite{BM}] $\dim(\geo \Ga)= cd_\Q(\Ga) -1$. 
\end{thm}

In particular, by Stallings--Swan--Dunwoody Theorem, $\Ga$ is {\em virtually free} (i.e. contains a free subgroup of finite index) if and only if $\geo \Ga$ is zero-dimensional, if and only if $\geo \Ga$ is totally disconnected, equivalently, it is either empty, or consists of two points or is homeomorphic to the Cantor set. 

One classifies 1-dimensional boundaries of hyperbolic groups as follows:

\begin{thm}
[Kapovich--Kleiner, \cite{KK}] Suppose that $\Ga$ is a hyperbolic group with connected 1-dimensional Gromov boundary. Then either $\geo \Ga$ is homeomorphic to $S^1$, or $\Ga$ splits as a finite graph of groups with virtually cyclic edge groups\footnote{and, hence, its Gromov boundary can be inductively described using boundaries of vertex groups}, or $\geo \Ga$ is homeomorphic to the Sierpinski carpet or the Menger curve.  
\end{thm}

\begin{example}
{\em Hyperbolic von Dyck groups} $D(p,q,r)$,
$$
D(p,q,r)= \<a, b, c| a^p=b^q=c^r=1, abc=1\>, p^{-1}+ q^{-1} + r^{-1} <1. 
$$
These are hyperbolic groups with Gromov boundary homeomorphic to $S^1$. Moreover, each $D(p,q,r)$ admits a unique (up to conjugation in $\Isom(\H^2))$ isometric conformal action on the hyperbolic plane.  
\end{example}

{\bf Representations of von Dyck groups to $PU(2,1)$.} Given an element $g\in G=PU(2,1)$ we let $\zeta(g)$ denote the codimension in $G$ of the centralizer of $g$ in $G$. In other words, $\zeta(g)$ is the local dimension near $g$ of the subvariety of elements of $G$ conjugate to $g$. Thus, $\zeta(g)\ge 2$ for every $g\in G$. Furthermore, if $g$ is an involution then $\zeta(g)=4$. 
The paper  \cite{Weil}  by Andre Weil describes  the local geometry of the {\em character variety} 
$$
Hom(D(p,q,r), G)//G
$$
as follows:

Suppose that $\rho: D(p,q,r)\to G$ is a {\em generic} representation, i.e. one whose image has  trivial centralizer in $G$. For instance, any representation whose image is discrete, nonelementary, not stabilizing a complex geodesic, will satisfy this condition. Then, near $[\rho]$, the real-algebraic variety  $Hom(D(p,q,r), G)//G$ is smooth  of dimension
$$
\zeta(\rho(a)) + \zeta(\rho(b)) + \zeta(\rho(c)) - 2 \dim (G)= \zeta(\rho(a)) + \zeta(\rho(b)) + \zeta(\rho(c)) -16. 
$$
Assuming that $p=2$, $\zeta(\rho(a))=4$, which implies that   
$$
\zeta(\rho(a)) + \zeta(\rho(b)) + \zeta(\rho(c)) -16 \le 4 +12-16=0. 
$$
Combined with an easy analysis  of non-generic representations, one obtains: 

\begin{proposition}\label{prop:rigid} 
If $p=2$ then $Hom(D(p,q,r), G)//G$ is zero-dimensional. 
\end{proposition}

\begin{example}
[Polygon-groups] \label{ex:poly}
Fix two natural numbers $p\ge 5$ and $q\ge 3$. Define the {\em polygon-group} $\Ga_{p,q}$ via presentation
$$
\<a_1,...,a_p| a_i^q=1, [a_i, a_{i+1}]=1, i=1,...,p\>, 
$$
where $i$ is taken mod $p$. Each $\Ga_{p,q}$ is hyperbolic with $\geo \Ga_{p,q}$ homeomorphic to the Menger curve. 
\end{example}

Every $\Ga_{p,q}$ admits a canonical {\em reflection representation} $\rho_{p,q}$ to $PU(2,1)$ constructed as follows: 

Pick a real hyperbolic plane $\H^2_\R\subset \H^2_\C$ and a regular right-angled $p$-gon $P=z_1...z_p$ in $\H^2_\R$. Let $C_i$ denote the complex geodesic through the edge $z_i z_{i+1}$ of $P$ ($i$ is taken mod $p$). For each $i$ let $g_i$ be the order $q$ complex reflection with the fixing $C_i$, with the rotation in the hyperplanes normal to $C_i$ through the angle 
$2\pi/q$.  Then $[g_i, g_{i+1}]=1$ and, hence, we obtain a representation
$$
\rho_{p,q}: \Ga_{p,q}\to PU(2,1). 
$$

\section{Appendix D: Orbifolds}

The notion of {\em orbifold} is a generalization of the notion of a
{\em manifold} which appears naturally in the context of properly
discontinuous non-free actions of groups on manifolds. Orbifolds were
first invented by Satake \cite{Satake} in 1950-s under the name of
V-manifolds, they were reinvented under the name of orbifolds by
Thurston in 1970's (see \cite{Thurston}) as a technical tool for proving his Hyperbolization
Theorem. We refer the reader to \cite{BMP} for a detailed treatment of orbifolds. 

\medskip
Before giving a formal definition we start with  basic examples of
orbifolds. Suppose that $M$ is a smooth connected manifold and $G$ is a discrete
group acting smoothly, faithfully\footnote{i.e. each nontrivial element of $G$ acts nontrivially}  and properly discontinuously on
$M$.  Then the quotient $\O= M/G$ is an {\em orbifold}, such orbifolds
are called {\em good}. The quotient $M/G$,  considered as a topological
space $X_\O$, is the {\em underlying space} of this orbifold.  If $S$ is
a set of points in $M$ where the action of $G$ is not free, then its
projection $\Si = S/G$ is the {\em singular locus} of the orbifold $\O$.

To be more concrete, consider 2-dimensional orbifolds. Suppose that 
$M= \H_\R^2$ and $G$ is a discrete subgroup of $PSL(2, \R)$. Then
the quotient $\O= \H^2/G$
is a Riemann surface $X_\O$ with a discrete collection of {\em cone points}
$z_j$ which form the singular locus $\Si$ of the orbifold $\O$. The projection
$p: \H^2 \to \O$ is the {\em universal cover} of the orbifold $\O$.
The Riemann surface $X_\O$  has a natural hyperbolic metric which is
singular in the discrete set 
$\Si$.  Metrically, the points $z_j$ are characterized by the property that the total angles around these
points are $2\pi/n_j$. The numbers $n_j$ are the  orders of cyclic subgroups
$G_{z_j}$ of $G$ which stabilize the points in $p^{-1}(z_j)$,
they are called {\em the local isotropy groups}.
The projection $p$ is a {\em ramified covering} 
from the point of view of Riemann surfaces.  From the point of view
of orbifolds this is an (orbi) {\em covering}. Thus, the singular
locus of the orbifold $\O$ consists of the points $z_j$ in $\Si$ equipped
with the extra data: The $PSL(2,\R)$-conjugacy classes of the local isotropy
groups $G_{z_j}$ (of course, each local isotropy group $G_{z_j}$ is determined
by the number $n_j$).

\medskip
We now discuss the general definition. A (smooth) $n$-dimensional 
orbifold $\O$ is a pair: A Hausdorff paracompact topological space $X$ 
(which is called the
{\em underlying space} of $\O$ and is denoted $X_\O$) and an {\em orbifold-atlas} $A$ on $X$.
The atlas $A$ consists of: 

\begin{itemize}
\item A collection of open sets $U_i \subset X$,which is closed under
taking finite intersections,   such that $X=\bigcup_i U_i$. 
\item A collection of open sets $\tilde{U}_i \subset\R^n$.  
\item A collection of finite groups of diffeomorphisms $\Ga_j$
acting on $\tilde{U}_i$ so that each nontrivial element of $\Ga_j$ acts nontrivially on each component of $\tilde{U}_j$.
\item  A collection of  homeomorphisms
$$ \phi_i : U_i \to \tilde{U}_i/\Ga_i.  $$
\end{itemize}
 We require the atlas $A$ to behave well
under inclusions.  Namely, if $U_i \subset U_j$, then
there is a smooth embedding $$ \tilde\phi_{ij} : \tilde{U}_i \to  \tilde{U}_j $$
and a monomorphism $f_{ij}: \Ga_i \to \Ga_j$ such that $\tilde\phi_{ij}$ is
$f_{ij}$-equivariant.

The open sets $U_j$ are the  {\em coordinate } neighborhoods of the points $x \in U_j$ 
and $\tilde{U}_j$ are their {\em covering coordinate} neighborhoods.

Similarly to orbifolds, one defines the class of {\em orbifolds with
boundary}; just instead of {\em open} 
sets $\tilde{U}_j \subset \R^n$   use open subsets in
$$
\R^{n}_+ \cup \R^{n-1}= \{ (x_1,...,x_n) :
x_n \ge 0\}  .
$$
The {\em boundary} of such orbifold consists of points
$x \in X_\O$ which correspond to $\R^{n-1}$ under the identification
$U_i \cong \tilde{U}_i/\Ga_i$.  As in the case of manifolds, the boundary
of each orbifold is an orbifold without boundary. By abusing the notation, 
we will call an {\em orbifold with boundary} simply an {\em orbifold}.
A compact orbifold without boundary is called
{\em closed}. 

To each point $x\in X$ we associate a germ of action
of a finite group of diffeomorphisms $\Ga_x$ at a fixed point $\tilde{x}$. If $\phi_j(x)$ is covered by a point
$\tilde{x}_j\in \tilde{U}_j$, then we have the isotropy group $\Ga_{j,x}$ of
$\tilde{x}_j$ in $\Ga_j$. Note that if  $U_i \subset U_j$, then the map
$ \tilde\phi_{ij} : \tilde{U}_i \to  \tilde{U}_j$ induces an isomorphism from the germ  
of the action of $\Ga_{j,x}$ at $\tilde{x}_j$ to
the germ  of the action of $\Ga_{i,x}$ at $\tilde{x}_i$. Thus, we let the  
germ $(\Ga_x, \tilde{x})$ be the equivariant diffeomorphism class of the germ
$(\Ga_{j,x}, \tilde{x}_j)$. The group $\Ga_x$ is called the {\em local
isotropy group}  of $\O$ at $x$. The set of points $x$ with nontrivial 
local isotropy group is called {\em the singular locus} of $\O$ and
is denoted by $\Si_\O$. 
Note that the singular locus is nowhere dense in $X_\O$.
An orbifold with empty singular locus is called {\em nonsingular}
or a {\em manifold}.

\medskip 
The main source of examples of orbifolds is:  

\begin{example}
Let $M$ a smooth connected $n$-manifold and $\Ga$ is a discrete
group acting smoothly and faithfully on $M$. Then $X=M/\Ga$ has a natural orbifold structure. The atlas $A$ on  $X$ is given as follows: Each $y\in M$ admits a coordinate neighborhood $\tilde U$ (identified with an open subset of $\R^n$) such that for every $g\in \Ga$ either $g\tilde U\cap \tilde U=\emptyset$ or $g\in \Ga_y$ (the stabilizer of $y$ in $\Ga$) 
and $g(\tilde U)=\tilde U$. Then let $\phi: \tilde{U}\to U= \phi(\tilde U)$ 
be the quotient map. One verifies that $A$ indeed satisfies axioms of an orbifold-atlas. The groups $\Ga_j$ in the definition of an atlas are just the stabilizers $\Ga_y$  as above. 
\end{example}

Since $\Ga_x$ acts smoothly near the fixed point $\tilde{x}$, the germ
$(\Ga_x, \tilde{x})$ is linearizable: We equip a neighborhood of $\tilde{x}$ with a $\Ga_x$-invariant
Riemannian metric; then the exponential map (with the origin at $\tilde{x}$)
conjugates the orthogonal action of $\Ga_x$ on the tangent space
$T_{\tilde{x}} \R^n$ to the germ of the action of $\Ga_x$ at $\tilde{x}$. 

\begin{defn}
A {\em Riemannian metric} $\rho$ on orbifold $O$ is the usual Riemannian
metric on $X_\O- \Si_\O$, such that after we lift $\rho$ to the local
covering coordinate neighborhoods $\tilde{U}_i$, it extends to a $\Ga_i$-invariant
Riemannian metric on the whole $\tilde{U}_i$.  
\end{defn}

The same definition applies to {\em complex structures}.  

\begin{exe}
Each orbifold $\O$ admits a Riemannian metric. Hint: use the partition of
unity argument similar to the manifold case.
\end{exe}

A {\em homeomorphism} (resp. {\em diffeomorphism}) between orbifolds
$\O, \O'$ is a homeomorphism
$h: X_\O \to X_{O'}$ such  that for all points $x \in \O, y= h(x)\in O'$, there
are coordinate neighborhoods $U_x \cong \tilde{U}_x/\Ga_x, V_y\cong \tilde{V}_y/\Ga_y$
such that $h$ lifts to an equivariant homeomorphism (resp.  diffeomorphism)
$$ \tilde{h}_{xy}: 
\tilde{U}_x \to \tilde{V}_y.  $$

Note that to describe a smooth orbifold $\O$ up to homeomorphism
it suffices to describe the topology of the pair $(X_\O, \Si_\O)$ and
the homeomorphic equivalence classes of the germs $(\Ga_x, \tilde{x})$
for the points $x\in \Si_\O$.

\begin{exe}
\label{segment}
Let $\O$ be a connected compact 1-dimensional orbifold without boundary which
is not a manifold. Then $\O$ is homeomorphic to the closed interval $[a,b]$
where $(\Ga_a, \tilde{a}), (\Ga_b, \tilde{b})$ are the germs $(\Z_2, 0)$ of the reflection
group $\Z_2$ acting isometrically on $\R$ near its fixed point $0\in \R$.  
\end{exe}

A {\em smooth} map between orbifolds $\O$ and $\O'$ is a continuous map
$$ g: \O \to O' $$
which can be (locally) lifted
to smooth equivariant maps between pairs of coordinate covering neighborhoods
$$ \tilde{g}_{ij}:
\tilde{U}_j \to \tilde{V}_i$$ 
Similarly we define {\em immersions} and {\em submersions} between orbifolds as smooth maps between orbifolds which locally
lift to immersions and submersions respectively.

% A {\em suborbifold}\index{suborbifold}  $R \subset O$ is an orbifold $R$ together with an {\em embedding} (i.e. an {\em injective immersion}) $\iota : R \hook O$. We say that a suborbifold $R$ is {\em properly embedded} in the orbifold $O$ if the embedding $\iota : R \hook O$ sends $\D R$ to $\D O$, $\iota^{-1}(\D O) \subset \D R$ and written in local covering coordinates the map $\iota$ is transversal to $\D O$.

\medskip Suppose that $\O', \O$ are orbifolds and $p: X_{\O'} \to X_{\O}$
is a continuous map. The map $p$ is called a {\em covering map}
between the orbifolds  $\O', \O$ if the following property is satisfied:

For each point $x \in  X_{\O}$ there exists a chart $U=\tilde{U}/\Ga_x$
such that for every component $V_i$ of $p^{-1}(U)$, the
restriction map $p: V_i \to U$ is a quotient map of an equivariant
diffeomorphism $h_i: \tilde{V}_i\to \tilde{U}$ (if $y_i= p^{-1}(x)\cap V_i$ then
$h_i$ conjugates the action of $\Ga_{y_i}$ on $\tilde{V}_i$ to the action of a
subgroup of $\Ga_x$ on $\tilde{U}$).

From now on we will assume that the orbifolds under consideration are connected.

\medskip
The {\em universal covering} $p:\tilde{\O}\to \O$ of an orbifold $\O$ is
the {\em initial object} in the category of orbifold coverings,
i.e. it is a covering such that
for any other covering $p': \O' \to \O$  
there exists a covering $\tilde{p}: \tilde{\O} \to \O'$ satisfying 
$p'\circ \tilde{p} = p$. If $p:\tilde{\O}\to \O$ is the universal covering
then  the orbifold $\tilde{\O}$ is called   the 
{\em universal covering orbifold} of $\O$.

The group $Deck(p)$ of {\em deck transformations} of an orbifold
covering $p: \O' \to \O$ is the group of self-diffeomorphisms $h: \O'
\to \O'$ such that $p \circ h = p$. A covering $p: \O' \to \O$ is called
{\em regular} if $\O'/Deck(p)= \O$. 

The {\em fundamental group} $\pi_1(\O)$ of the orbifold $\O$ is the
group of deck transformations of its universal covering. Then $\O=
\tilde{\O}/\pi_1(\O)$. An alternative definition of the fundamental
group based on homotopy-classes of loops in $O$ see in
\cite[Chapter 13]{Ratcliffe}. 

%The following theorem was proven by Thurston in \cite[Chapter 13]{Thurston(1978-81)}.

\begin{thm}
Each orbifold has a universal covering.
\end{thm}

\begin{defn}
 An orbifold $\O$ is called {\em good} if its
universal covering is a manifold. Orbifolds which are not good are
called {\em bad}. An orbifold is called {\em very good}
if is admits a finite-sheeted manifold-covering space.  
\end{defn}

\begin{example}
Let $\O=M_\Ga$ be an $n$-dimensional complex hyperbolic orbifold. Then $\Ga=\pi_1(\O)$ and $\O$ is a good orbifold: Its universal covering space is $\H^n_\C$.  If $\Ga$ is finitely generated 
then, according to Selberg's Lemma, the orbifold $\O$ is very good. 
\end{example}

{\bf Orbifold bundles.} Instead of defining orbifold bundles in full generality, I will define  these only in the case of compact fibers andf connected base, since this will suffice for our purposes:

\begin{defn}
A smooth orbi-bundle with compact fibers and connected base is a proper submersion $f: \O\to \BB$ between orbifolds. 
Fibers of $f$ are preimages of points under $f$. 
\end{defn} 

Note that two different fibers need not be isomorphic to each other, but one can prove that they are {\em commensurable} in the sense that they have a common finite-sheeted orbi-covering.

\section{Appendix E: Ends of spaces}

Let $Z$ be a locally path-connected, locally compact, Hausdorff topological space. The set of {\em ends} of $Z$ can be defined as follows (see e.g. \cite{DK} for details). 

Consider an exhaustion $(K_i)$ of $Z$ by an increasing sequence of   compact subsets:
$$K_{i}\subset K_{j}, \quad \textup{ whenever } i\leq j,$$
and 
$$\bigcup_{i\in \mathbb{N}} K_{i}=Z.$$
Set $K_i^c:= Z\setminus K_i$. The ends of $Z$ are equivalence classes of decreasing sequences of connected components 
$(C_i)$ of ${K_i}^c$: 
$$C_{1}\supset C_{2}\supset C_{3}\supset \cdots $$
Two sequences $(C_i), (C'_j)$ of components of $({K_i}^c), ({K'_j}^c)$ are said to be equivalent if each $C_{i}$ contains some $C'_{j}$ and vice-versa. Then the equivalence class of a sequence $(C_i)$ is an {\em end} $e$ of $Z$. Each $C_i$ and its closure is called a 
 \emph{neighborhood} of $e$ in $Z$. The set of ends of $Z$ is denoted $Ends(Z)$. 
 An end $e$ is called {\em isolated} if it admits a closed 1-ended neighborhood $C$; such a neighborhood is called {\em isolating}.  
 
 An alternative view-point on the neighborhoods of ends is that 
 there is a natural topology on the union $\hat{Z}=Z\cup Ends(Z)$  which is a compactification of $Z$ and the neighborhoods $C$ of ends $e$ above are intersections of $Z$ with neighborhoods of  $e$ in $\hat{Z}$. Then an end $e$ is isolated if and only if it is an isolated point of $\hat{Z}$.  
 A closed neighborhood $C$ of $e$ in $Z$ is isolating if and only if $C\cup \{e\}$ is closed in  $\hat{Z}$. 
 
From this definition it is not immediate that the notion of ends is independent on the choice of an exhausting sequence $(K_i)$ of compact subsets. The true, but less intuitive, definition of 
$Ends(Z)$ is by considering the poset (ordered by the inclusion) of all compact subsets $K\Subset Z$. This poset defines the inverse system of sets  
$$
\{\pi_0(K^c, x): K\Subset Z\},
$$
where the inclusion $K\subset K'$ induces the map 
$$
\pi_0(Z-K', x')\to \pi_0(Z-K, x'), 
$$
with $x'\in Z-K'\subset Z-K$. Taking the inverse limit of this system of  sets yields $Ends(Z)$ which is, manifestly, a topological invariant. Furthermore, it is an invariant of the proper homotopy type of $Z$. 
 
In this lectures, I  adopt the {\em analyst's viewpoint} on ends of manifolds and conflate isolated ends and their isolating neighborhoods.

\section{Appendix F: Generalities on function theory on complex manifolds} \label{sec:CA} 

For a complex manifold $M$ let $\O_M$ denote the ring of holomorphic functions on $M$. 
By a {\em complex manifold with boundary} $M$ I mean a smooth manifold with (possibly empty) boundary $\partial M$,  
such that the interior, $\inte(M)$, of the manifold $M$, is equipped with a complex structure, and there exists a smooth embedding $h: M\to X$ to an equidimensional complex manifold $X$,  biholomorphic on $\inte(M)$.  
A holomorphic function on $M$ is a smooth function which admits a holomorphic extension to a neighborhood of $M$ in $X$. 

Suppose that $X$ is a complex manifold and $Y\subset X$ is a codimension $0$ smooth submanifold with boundary in $X$. The submanifold 
$Y$ is said to be {\em strictly Levi-convex} if  every boundary point of $Y$ 
admits a neighborhood $U$ in $X$ such that the submanifold with boundary $Y\cap U$ can be written  as
$$
\{\phi\le 0\},
$$
for some smooth submersion $\phi: U\to \R$ satisfying 
$$ 
Hess(\phi) >0.   
$$
Here $Hess(\phi)$ is the holomorphic Hessian:
$$
\left( \frac{\partial^2 \phi}{\partial \bar{z}_i \partial z_j}\right).  
$$
(Positivity of  the Hessian is independent of the local holomorphic coordinates.) 

\begin{example}
If $X=\C^n, Y=\{z\in \C^n: |z| \le 1\}$, then $Y$ is strictly Levi-convex in $X$: The complex Hessian of the function $\phi(z)= |z|^2= z\cdot \bar{z}$ is the identity matrix. 
\end{example}

%Hormander, Thm 2.6.10

\begin{defn}
A {\em strongly pseudoconvex manifold} $M$ is a complex manifold with boundary 
which admits a strictly Levi-convex holomorphic embedding in an equidimensional complex manifold. 
\end{defn}

Suppose, in addition, that $M$ is compact and $h: M\to X$ is a holomorphic embedding with strictly Levi-convex image $Y$. 
Then there exists a strictly Levi-convex submanifold $Y'\subset X$ such that $Y\subset \inte(Y')$. Accordingly, $M$ can be biholomorphically embedded in the interior of a compact strongly pseudoconvex manifold $M'$.

\begin{defn}
An  complex manifold $Z$ is called {\em holomorphically convex} if for every discrete closed subset $A\subset Z$ there exists a holomorphic function $Z\to \C$ which is proper on $A$. 
\end{defn}

%... Discuss equivalent definitions...
Alternatively,\footnote{and this is the standard definition} one can define holomorphically convex manifolds as follows: For a compact $K$ 
in a complex manifold $M$, the {\em holomorphic convex hull} $\hat{K}_M$ of $K$ in $M$ is 
$$
\hat{K}_M= \{z\in M: |f(z)|\le \sup_{w\in K} |f(w)|, \forall f\in \O_M\}.  
$$
Then $M$ is holomorphically convex iff for every compact $K\subset M$, the hull $\hat{K}_M$ is also compact. 

\begin{defn}
A complex manifold is called {\em Stein} if it admits a proper holomorphic embedding in $\C^N$ for some $N$. 
\end{defn}

Equivalently, $M$ is Stein iff it is holomorphically convex and any two distinct points $z, w\in M$ can be separated by a holomorphic function, i.e. there exists $f\in \O_M$ such that $f(z)\ne f(w)$. Yet another equivalent definition is: A complex manifold $M$ is Stein if and only if it is strongly pseudoconvex, i.e. it admits an exhaustion by codimension 0 
strongly pseudoconvex complex submanifolds with boundary. 

In particular:

\begin{thm}
%[Grauert] 
The interior of every compact strongly pseudoconvex manifold $Z$ is holomorphically convex. 
\end{thm}

Therefore, by holomorphically embedding such (connected manifold) $Z$ in the interior of another compact strongly pseudoconvex manifold $Z'$ and applying Grauert's theorem to $Z'$, it follows that $Z$ admits nonconstant holomorphic functions.

Kohn and Rossi in \cite{KR} proved a certain extension theorem for CR functions defined on the boundary of a complex manifold 
 to holomorphic functions on the entire manifold. I will state only a weak form of their result which will suffice for our purposes. 

\begin{thm}
[Kohn--Rossi] \label{thm:Kohn-Rossi} 
Suppose that $M$ is a compact strongly pseudoconvex complex manifold of dimension $>1$ which admits at least one 
nonconstant holomorphic function. Then every holomorphic function on $\partial M$ extends to a holomorphic function on the entire $M$. 
\end{thm}

As one of the corollaries of this theorem (Corollary 7.3 of \cite{KR}), it follows  that if such $M$ is connected then $\partial M$ is also connected. (If $\partial M$ is disconnected, then one can take a nonconstant locally constant function defined near $\partial M$: Such a function cannot have a holomorphic extension to $M$.) 

\begin{rem}
If $M$ is K\"ahler, then Theorem \ref{thm:Kohn-Rossi} also holds without the assumption on the existence of nonconstant holomorphic functions, see Proposition 4.4 in \cite{NR}. 
\end{rem}

\begin{thm}
[Rossi, \cite{R2}, Corollary on page 20] \label{thm:rossi} 
Suppose that $M$ is a compact strongly pseudoconvex  complex manifold. Then $\inte(M)$ admits a proper surjective holomorphic map 
to a Stein space. In particular, if $\inte(M)$ contains no compact complex subvarieties of positive dimension, then $\inte(M)$ is Stein. 
\end{thm} 

I will not define Stein spaces here (strictly speaking, they are not needed for the purpose of these notes), I refer to \cite{GR} for various equivalent definitions. 

\medskip 
{\bf Topology of Stein manifolds and spaces.} 
Every complex $n$-dimensional Stein space is homotopy-equivalent to an $n$-dimensional CW complex, see \cite{Hamm1, Hamm2}. More precisely (see Theorem 1.1* on page 153 of \cite{GM}):

\begin{thm}
Let $M$ be a  $n$-dimensional complex manifold which admits a proper holomorphic map $M\to \C^N$ 
with fibers of positive codimension. Then $M$ is homotopy-equivalent to an $n$-dimensional CW complex. 
\end{thm}

\begin{cor}\label{cor:CW}
Suppose that $M$ is a connected compact strongly pseudoconvex  complex $n$-manifold with nonempty boundary. Then   $M$ is homotopy-equivalent to a CW complex of dimension $2n-2$. 
\end{cor}

\section{Appendix G (by Mohan Ramachandran): Proof of Theorem \ref{thm:mohan}} \label{sec:mohan}

\begin{prop}
Let $X$ be a complex manifold of dimension $\ge 2$ and let 
$M\subset X$ be a domain with compact nonempty smooth strongly pseudoconvex 
boundary. Then every pluriharmonic function on 
$M$ which vanishes at $\partial M$,  vanishes identically. 
\end{prop}
\proof The proof mostly follows that of Proposition 4.4 in \cite{NR}. 
 Suppose that $M=\{x\in X: \varphi(x)<0\}$ for some smooth function $\varphi$, which is strictly plurisubharmonic on a neighborhood of $\partial M$ and such that 
 there exists $\eps<0$ such that $\varphi^{-1}([\eps, 0])$ is compact and $\varphi|_{\partial M}=0$.  
 Let $\beta: M\to \R$ be a pluriharmonic function  which vanishes at $\partial M$.
Fix $a\in (\eps, 0)$, such that $\varphi$ is strictly plurisubharmonic  on $V=\{x\in M: \varphi(x)>a\}$. 
%If $\rho$ is the real or imaginary part of $\beta$ and $\rho$ does not vanish 
If $\beta$ does not vanish identically on a neighborhood of $\partial M$, we let $b\in \beta(V)$ 
denote a regular value of $\beta$.  Thus, $\beta^{-1}(b)$ is disjoint from $\partial M$. 
Since $\varphi^{-1}([\eps, 0])$ is compact, the restriction 
of $\varphi$ to $\beta^{-1}(b)$ has a maximum at some point $x_0\in V\cap \beta^{-1}(b)$. The holomorphic $1$-form $\partial \beta$ 
determines a (singular) holomorphic foliation on $M$. 
Consider the leaf $L$ through $x_0$ of this holomorphic foliation: This leaf is contained in 
$\beta^{-1}(b)$ and, hence, the restriction $\varphi|_L$  has a maximum at $x_0$ contradicting strict plurisubharmonicity of $\varphi$. Therefore, $\beta$ is identically zero  near $\partial M$ and, hence, is identically zero. \qed 

\medskip 
The next proposition is proven in \cite[Theorem 2.6]{NR1}: 

\begin{prop}
Suppose now that $M$ has a complete K\"ahler metric of bounded geometry\footnote{i.e. its sectional curvature lies in a finite interval and its injectivity radius is bounded from below}, $\partial M$ is connected and $M$ has at least two ends. Then $M$ admits 
a nonconstant pluriharmonic function $\beta: M\to \R$ which converges to zero 
at $\partial M$. 
\end{prop}

\medskip 
By combining the two propositions, we conclude: 

\begin{cor}
Suppose that $M$ is a complex manifold of dimension $\ge 2$, which admits a holomorphic embedding 
as a domain with compact nonempty smooth strongly pseudoconvex 
boundary and which admits a complete K\"ahler metric of bounded geometry. Then $M$ is 1-ended. 
\end{cor}

We can now conclude the proof of  Theorem \ref{thm:mohan}: Let $M=M_\Ga$ be a complex hyperbolic manifold of dimension $\ge 2$ and of injectivity radius bounded below. Suppose that $E_0\subset M$ is a convex end. Let $S_0\subset \partial \ol{M}$ be the component corresponding to the end $E_0$. Consider the complex manifold $Y=\check{\Omega}_\Ga/\Ga$. Remove from $Y$ all the components of $Y-M$ which are disjoint from $S_0$ and call the result $X$. Then $M$ embeds in $X$ as a domain with 
nonempty smooth strongly pseudoconvex boundary, namely, $S_0$.  Then, by the corollary, $M$ is 1-ended. 
 \qed


\newpage 
\begin{thebibliography}{ABC} 

\bibitem{AGG} 
A. Anan'in, C. H. Grossi and  N. Gusevskii, {\em Complex hyperbolic structures on disc bundles over surfaces}, International Mathematics Research Notices, Vol. 19 (2011) pp. 4295--4375.

\bibitem{Ancona}
A. Ancona, {\em Convexity at infinity and Brownian motion on manifolds with unbounded negative curvature}, Rev. Mat. Iberoamericana, Vol. 10 (1994), pp. 189--220.


\bibitem{Anderson} 
M. Anderson, {\em The Dirichlet problem at infinity for manifolds of negative curvature},
J. Differential Geometry, Vol. 18 (1983), pp. 701--721


 \bibitem{Ballmann}
  W. Ballmann, ``Lectures on Spaces of Nonpositive Curvature.'' With an appendix by Misha Brin. 
DMV Seminar, vol.\ 25, Birkh\"auser Verlag, Basel, 1995. 

  
  \bibitem{BGS} W. Ballmann, M. Gromov and V. Schroeder, ``Manifolds of nonpositive curvature.'' Progr. Math. 61, Birkh\"auser, Boston, 1985. 
   



\bibitem{BHH} 
G. Barthel, F. Hirzebruch and T. H\"{o}fer, 
``Geradenkonfigurationen und Algebraische Flächen.''  %[Line Configurations and Algebraic Surfaces]
Aspects of Mathematics, D4. Friedr. Vieweg \& Sohn, Braunschweig, 1987.

\bibitem{BK} 
I. Belegradek and V. Kapovitch, {\em Classification of negatively pinched manifolds with amenable fundamental groups}, 
Acta Math., Vol. 196 (2006) pp. 229--260.




\bibitem{BM}
M. Bestvina and G. Mess, {\em The boundary of negatively curved groups}, J. Amer. Math. Soc. Vol. 4 (1991) pp. 469--481.

\bibitem{BMP}
M. Boileau, S. Maillot and J. Porti, ``Three-dimensional Orbifolds and Their Geometric Structures.'' Panoramas and Syntheses 15. Soci\'et\'e Math\'ematique de France, 2003.


\bibitem{Borbely}
A. Borb\'ely,  {\em Some results on the convex hull of finitely many convex sets}, 
Proc. Amer. Math. Soc. Vol. 126 (1998), no. 5, pp. 1515--1525.


\bibitem{Bourdon} 
M. Bourdon, Sur la dimension de Hausdorff de l'ensemble limite d'une famille
de sous-groupes convexes co-compacts, C. R. Acad. Sci. Paris S\'er. I Math., Vol. 325
(1997), pp. 1097--1100.

\bibitem{Bowditch0} B.H. Bowditch, {\em Some results on the geometry of convex hulls in manifolds of pinched
negative curvature}, Comment. Math. Helv., Vol. 69 (1994), pp. 49--81.

\bibitem{Bowditch}
B. H. Bowditch, {\em Geometrical finiteness with variable negative curvature}, Duke Math. Journal, Vol. 77 (1995), no. 1, pp. 229--274. 

\bibitem{Bowditch1999}
B. H. Bowditch, 
{\em Convergence groups and configuration spaces}. 
In: ``Geometric group theory down under'' (Canberra, 1996), 
pp.\ 23--54, de Gruyter, Berlin, 1999. 

\bibitem{Bowen}
L. Bowen,  
{\em Cheeger constants and $L^2$-Betti numbers},   
Duke Math. Journal, Vol. 164 (2015), no. 3, pp. 569--615.

\bibitem{BI}
M. Burger and A. Iozzi,  {\em A measurable Cartan theorem and applications to deformation rigidity in complex hyperbolic geometry}, 
Pure Appl. Math. Q. Vol. 4 (2008), no. 1, Special Issue: In honor of Grigory Margulis. Part 2, pp. 181--202.


\bibitem{BS} 
D. Burns and S. Shnider, {\em Spherical hypersurfaces in complex manifolds}, Invent. Math., Vol. 33
(1976) pp. 223--246. 



\bibitem{CNS}
A. Cano, J. Navarrete and J. Seade, ``Complex Kleinian Groups,'' Progress in Mathematics, Vol. 
303, Springer Verlag, 2013. 

\bibitem{CS}
A. Cano and J. Seade, {\em On the equicontinuity region of discrete subgroups of $PU(1,n)$}, 
Journal of Geometric Analysis Vol. 20 (2), pp. 291--305. 

\bibitem{CMP}
J. Carlson, S. M\"uller-Stach and C. Peters, 
``Period Mappings and Period Domains.'' Cambridge Studies in Advanced Mathematics, 2003. 

\bibitem{Coo}
M. Coornaert, {\em Mesures de Patterson-Sullivan sur le bord d'un espace hyperbolique au sens de Gromov}, Pacific Journal of Mathematics, Vol. 159 (1993), pp. 241--270. 

\bibitem{Corlette0}
K. Corlette, {\em Flat $G$-bundles with canonical metrics}, 
J. Differential Geom. Vol. 28 (1988), no. 3, pp. 361--382.

\bibitem{Corlette}
K. Corlette, {\em Hausdorff dimensions of limit sets I,} Inventiones Math. Vol. 102 
(1990), pp. 521--541.


\bibitem{Corlette-NA}
K. Corlette, {\em Archimedean superrigidity and hyperbolic geometry}, Ann. Math. (2) Vol. 135 (1992), pp. 165--182. 


\bibitem{CO}
K. Corlette and A. Iozzi, {\em  Limit sets of discrete groups of isometries of exotic hyperbolic spaces}, Trans. Amer. Math. Soc. Vol. 351 (1999), no. 4, pp. 1507--1530.

\bibitem{CHL} 
W. Couwenberg, G. Heckman and E. Looijenga, {\em Geometric structures on the complement of a projective arrangement}, 
Publ. Math. Inst. Hautes \'Etudes Sci., No. 101 (2005), pp. 69--161.


\bibitem{DOP}
F. Dal'bo, J.-P. Otal and  M.  Peign\'e, {S\'{e}ries de {P}oincar\'{e} des groupes g\'{e}om\'{e}triquement finis}, 
Israel J. Math., Vol. 118 (2000), pp. 109--124.

\bibitem{DSU}
T. Das, D. Simmons and M. Urba\'nski, ``Geometry and Dynamics in Gromov Hyperbolic Metric Spaces. With an Emphasis on Non-proper Settings.'' Mathematical Surveys and Monographs, 218. American Mathematical Society, Providence, RI, 2017. 

\bibitem{Delzant}
T. Delzant,  {\em L'invariant de Bieri--Neumann--Strebel des groupes fondamentaux des vari\'et\'es K\"ahl\'eriennes},  
Math. Ann. Vol. 348 (2010), no. 1, pp. 119--125.


\bibitem{DH} M. Desgroseilliers and F. Haglund, {\em On some convex cocompact groups in real hyperbolic space}, 
Geom. Topol. Vol. 17  (2013), pp. 2431--2484.

\bibitem{DM}
P. Deligne, G. D. Mostow, ``Commensurabilities Among Lattices in $PU(1,n)$.'' Annals of Mathematics Studies, vol. 132. Princeton University Press, Princeton, 1993. 

\bibitem{Deraux}
M. Deraux, {\em A new non-arithmetic lattice in $PU(3,1)$}, Algebr. Geom. Topol., to appear. 


\bibitem{DeK}
S. Dey, M. Kapovich, in preparation. 

\bibitem{DPP}
M. Deraux,  J. Parker, and J. Paupert, 
{\em New non-arithmetic complex hyperbolic lattices}, 
Invent. Math. Vol. 203 (2016), no. 3, pp. 681--771. 


\bibitem{Epstein}
D. B. A. Epstein, {\em Complex hyperbolic geometry}, In: ``Analytical and geometric aspects of hyperbolic space (Coventry/Durham, 1984),''  pp. 93--111, Cambridge Univ. Press, 1987. 

\bibitem{Goldman}
W. Goldman, {\em Complex hyperbolic Kleinian groups.} In: ``Complex Geometry (Osaka, 1990), pp. 31--52.'' Lecture Notes in Pure and Appl. Math., Vol. 143, Dekker, New York, 1993.
 
 
\bibitem{Gol}
W. Goldman, ``Complex Hyperbolic Geometry,'' Oxford Mathematical Monographs, 1999. 

\bibitem{GKL}
W. Goldman, M. Kapovich and B. Leeb, Complex hyperbolic surfaces homotopy-equivalent to a Riemann surface, Comm. in Analysis and Geom., Vol. 9 (2001), pp. 61--96. 

\bibitem{GoMi}
W. Goldman and J. J. Millson, {\em Local rigidity of discrete groups acting on complex hyperbolic space}, Invent. Math. Vol. 88 (1987), no. 3, pp. 495--520.

\bibitem{GM}
M. Goresky and R. MacPherson, ``Stratified Morse Theory.'' 
Ergebnisse der Mathematik und ihrer Grenzgebiete, Vol. 14. Springer-Verlag, Berlin, 1988.


\bibitem{Gran}
J. Granier, ``Groupes discrets en g\'eom\'etrie hyperbolique -- Aspects effectifs.'' PhD thesis, Universit\'e de Fribourg, 2015. 

\bibitem{Grau}
H. Grauert, {\em On Levi's problem and the imbedding of real-analytic manifolds},  
Ann. of Math. (2) Vol. 68 (1958) pp. 460--472. 

\bibitem{GR}
H. Grauert and  R. Remmert, ``Stein Spaces,'' Grundlehren der Mathematischen Wissenschaften, 236, Berlin-New York: 
Springer-Verlag, 1979. 

\bibitem{GPS}
M. Gromov and I. Piatetski-Shapiro, {\em Nonarithmetic groups in Lobachevsky spaces}, Inst. Hautes \'Etudes Sci. Publ. Math., No. 66 (1988), pp. 93--103. 

\bibitem{GS}
M. Gromov and R. Schoen, {\em Harmonic maps into singular spaces and $p$-adic superrigidity for lattices in groups of rank one}, Inst. Hautes \'Etudes Sci. Publ. Math., No. 76 (1992), pp. 165--246 

\bibitem{DK}
C. Drutu and M. Kapovich, ``Geometric Group Theory." AMS  Colloquium Publications,  2018.

\bibitem{Hamm1}
H. Hamm,  {\em Zur Homotopietyp Steinscher R\"aume},  
J. Reine Angew. Math. Vol. 338 (1983), pp. 121--135.

\bibitem{Hamm2}
H. Hamm,  {\em Zum Homotopietyp q-vollst\"andiger R\"aume}, J. Reine Angew. Math. Vol. 364 (1986), pp. 1--9.
 
 
 
\bibitem{HK} 
J. A. Hillman and D. H. Kochloukova, {\em Finiteness conditions and $PD(r)$-group covers of
$PD(n)$-complexes}, Math. Z., Vol. 256 (2007), pp. 45--56.


\bibitem{Hirz} F. Hirzebruch, {\em Arrangements of lines and algebraic surfaces}. In: ``Arithmetic and geometry,'' Vol. II, pp. 113--140,
Progr. Math., 36, Birkh\"auser, Boston, Mass., 1983. 

\bibitem{Holz} R.-P. Holzapfel, ``Ball and Surface Arithmetics.''
Aspects of Mathematics, E29. Friedr. Vieweg \& Sohn, Braunschweig, 1998.

\bibitem{Kapovich1998}
M. Kapovich, {\em On normal subgroups in the fundamental groups of complex surfaces}, Preprint,
math.GT/9808085, 1998. 

\bibitem{K2000}
M. Kapovich, ``Hyperbolic Manifolds and Discrete Groups,'' Birkhauser's series ``Progress in Mathematics,'' 2000. Reprinted in 2009. 

\bibitem{Kapovich2005}
M. Kapovich, {\em Representations of polygons of finite groups}, Geometry and Topology, Vol. 9 (2005) pp. 1915--1951. 

\bibitem{Kapovich2008}
M. Kapovich,  {\em Kleinian groups in higher dimensions}, In: ``Geometry and dynamics of groups and spaces,'' pp. 487--564, Progr. Math., Vol. 265, Birkh\"auser, Basel, 2008.
 
 \bibitem{K2018}
M. Kapovich,  {\em A note on Selberg's lemma and negatively pinched Hadamard manifolds}, ArXiv, 1808.01602, 2018.

\bibitem{proper}
M. Kapovich,  {\em A note on properly discontinuous actions}, Preprint, 2017. 


\bibitem{KK}
M. Kapovich and B. Kleiner, {\em Hyperbolic groups with low-dimensional boundary}, Ann. Sci. \'Ecole Norm. Sup. (4) Vol. 33 (2000), no. 5, pp. 647--669. 

\bibitem{KLP}
M. Kapovich, B. Leeb and J. Porti, {\em Dynamics on flag manifolds: domains of proper discontinuity and cocompactness},  
Geometry and Topology, Vol. 22 (2017) pp. 157--234.

\bibitem{KL}
M. Kapovich and B. Liu, {\em Geometric finiteness in negatively pinched Hadamard manifolds}, Ann. Acad. Sci. Fennicae, Vol. 44 (2019) 2, pp. 841--875.


\bibitem{Kazhdan}
D. Kazhdan, {\em Some applications of the Weil representation}, J. Analyse Mat., Vol. 32 (1977), pp. 235--248.


\bibitem{KR} J. J. Kohn and H. Rossi, {\em On the extension of holomorphic functions from the boundary of a complex manifold},  
Ann. of Math. (2) Vol. 81 (1965) pp. 451--472.

\bibitem{KM}
V. Koziarz and  N. Mok, 
{\em Nonexistence of holomorphic submersions between complex unit balls equivariant with respect to a lattice and their generalizations},  Amer. J. Math. Vol. 132 (2010), no. 5, pp. 1347--1363.

\bibitem{Krantz}
S. Krantz,  ``Function Theory of Several Complex Variables.'' 
Reprint of the 1992 edition. AMS Chelsea Publishing, Providence, RI, 2001. 

\bibitem{Ledrappier}
F. Ledrappier, {\em Structure au bord des vari\'{e}t\'{e}s {\`a} courbure n\'{e}gative}, 
In {\em ``S\'{e}minaire de {T}h\'{e}orie {S}pectrale et
  {G}\'{e}om\'{e}trie, {N}o. 13, {A}nn\'{e}e 1994--1995''}, Vol.~13 of {\em
  S\'{e}min. Th\'{e}or. Spectr. G\'{e}om.}, pp. 97--122. Univ. Grenoble I,
  Saint-Martin-d'H{\`e}res, 1995.


\bibitem{Leuz}
E. Leuzinger, {\em Critical exponents of discrete groups and $L^2$-spectrum}, Proc. Amer. Math. Soc., Vol. 132 (2004), no. 3, pp. 919--927.

\bibitem{Liu}
K. Liu, {\em Geometric height inequalities}, Math. Res. Lett., Vol. 3 (1996), pp. 693--702.


\bibitem{Livne}
R. Livne, ``On certain covers of the universal elliptic curve.''  
Ph. {D}. {T}hesis, {H}arvard {U}niversity, 1981. 

\bibitem{Margulis}
G. A. Margulis, {\em Discrete groups of motions of manifolds of nonpositive curvature.} In: ``Proceedings of the International Congress of Mathematicians (Vancouver, B.C., 1974),'' 
Vol. 2, pp. 21--34. Canadian Mathematical Congress, Montreal, 1975. 

\bibitem{Margulis-book}
G. A. Margulis, ``Discrete Subgroups of Semisimple Lie Groups.'' Ergebnisse der Mathematik und ihrer Grenzgebiete, Vol. 17. Springer, Berlin, 1991. 



\bibitem{McReynolds}
 D. B. McReynolds, {\em Arithmetic lattices in $SU(n,1)$}. Preprint, 2005. 

\bibitem{Mostow1980}
G. D. Mostow, {\em On a remarkable class of polyhedra in complex hyperbolic space}, Pac. J. Math. Vol. 86 (1980), pp. 171--276. 

\bibitem{NR1} 
T. Napier and M. Ramachandran, {\em Structure theorems for complete K\"ahler manifolds and applications to Lefschetz type theorems}, GAFA, Vol. 5 (1995) pp. 809--851. 



\bibitem{NR} T. Napier and M. Ramachandran, {\em The $L^2$ $\bar{\partial}$-method, weak Lefschetz theorems, and the topology of K\"ahler manifolds}, J. Amer. Math. Soc. Vol. 11 (1998), no. 2, pp. 375--396. 

\bibitem{Parker}
J. Parker, ``Hyperbolic Spaces.'' Jyv\"askyl\"a Lectures in Mathematics, Vol. 2, 2008. 

\bibitem{Parker2}
J. Parker,  {\em Complex hyperbolic lattices}. In: ``Discrete groups and geometric structures,'' pp. 1--42,
Contemp. Math., Vol. 501, Amer. Math. Soc., Providence, RI, 2009. 

\bibitem{Paupert}
J. Paupert, {\em Introduction to Hyperbolic Geometry}, Arizona State University Lecture Notes, 2016. 
 
\bibitem{Prokhorov}
M. Prokhorov, {\em Absence of discrete groups of reflections with a noncompact fundamental polyhedron of finite volume in a Lobachevskii space of high dimension}, 
Izv. Akad. Nauk SSSR Ser. Mat., Vol. 50 (1986), no. 2, pp. 413--424.

\bibitem{Ratcliffe}
J. Ratcliffe, ``Foundations of Hyperbolic Manifolds.'' Second edition. Graduate Texts in Mathematics, Vol. 149. Springer, New York, 2006.
 
\bibitem{RT}
T. Roblin and S. Tapie, {\em Critical exponent and bottom of the spectrum in pinched negative curvature}, Math. Res. Lett., Vol. 22 (2015), no. 3, pp. 929--944.
 
 
\bibitem{Rogawski}
J. Rogawski, ``Automorphic Representations of Unitary Groups in Three Variables.'' Vol. 123 of
Annals of Mathematics Studies, Princeton University Press, Princeton, NJ, 1990.


\bibitem{R1} H. Rossi, {\em Attaching analytic spaces to an analytic space along a pseudoconcave boundary}, In: ``Proc. Conf. Complex Analysis (Minneapolis, 1964),'' Springer Verlag, Berlin, 1965, pp. 242--256.  



\bibitem{R2} H. Rossi, 
{\em Strongly pseudoconvex manifolds},  In: ``Lectures in Modern Analysis and Applications, I,'' 
pp. 10--29. Springer Verlag, Berlin,  1969.  

\bibitem{Satake}
I. Satake, {\em On a generalization of the notion of manifold,} Proc. Nat. Acad. Sci. U.S.A., Vol. 42 (1956) pp. 359--363.

\bibitem{Shv}
O. Shvartsman,  {\em An example of a nonarithmetic discrete group in the complex ball.} In: ``Lie groups, their discrete subgroups, and invariant theory,'' pp. 191--196,
Adv. Soviet Math., Vol. 8, Amer. Math. Soc., Providence, RI, 1992. 


\bibitem{Sch}
R. Schwartz, ``Spherical CR Geometry and Dehn Surgery.'' Annals of Mathematics Studies, Vol. 165. Princeton University Press, Princeton, NJ, 2007. 

\bibitem{Simpson}
C. Simpson, {\em Higgs bundles and local systems}, Inst. Hautes \'Etudes Sci. Publ. Math. No. 75 (1992), pp. 5--95. 

\bibitem{Thurston}
W. Thurston, ``The Geometry and Topology of Three-Manifolds.'' Princeton University Lecture Notes. 1978--1981. Chapter 13. 

\bibitem{Th}
W. Thurston, {\em Shapes of polyhedra and triangulations of the sphere}, Geom. Topol. Monogr. Vol. 1 (1998), pp. 511--549. 

\bibitem{Toledo}
D. Toledo, {\em Harmonic maps from surfaces to certain K\"ahler manifolds}, Mathematica Scandinavica, Vol. 45 (1979)
pp. 13--26. 

\bibitem{Tukia1994}
P. Tukia, 
{\em Convergence groups and Gromov-hyperbolic metric spaces}, 
New Zealand Journal of Math.,  Vol. {23} (1994), pp. 157--187. 


\bibitem{Tukia1998}
P. Tukia, 
{\em Conical limit points and uniform convergence groups}, 
J.\ Reine Angew.\ Math.,  Vol. {501} (1998), pp. 71--98. 


\bibitem{Vinberg}
E. Vinberg,  {\em Absence of crystallographic groups of reflections in Lobachevskii spaces of large dimension}, Trudy Moskov. Mat. Obshch. Vol. 47 (1984), pp. 68--102.

\bibitem{Wallach} 
N. Wallach, {\em Square integrable automorphic forms and cohomology of arithmetic quotients of
$SU(p, q)$}, Math. Ann., Vol. 266 (1984), pp. 261--278.


\bibitem{Weil}
A. Weil, {\em Remarks on the cohomology of groups}, Annals of Mathematics (2), Vol. 80 (1964), pp. 149--157. 


\bibitem{Yau}
S.-T. Yau, {\em Calabi's conjecture and some new results in algebraic geometry}, 
Proc. Nat. Acad. Sci. U.S.A. Vol. 74 (1977), no. 5, pp. 1798--1799.

\bibitem{Yeung} 
S.-K. Yeung, {\em Virtual first Betti number and integrality of compact complex two-ball quotients},
Int. Math. Res. Not., Vol. 38 (2004), pp. 1967--1988.

\bibitem{Yue}
C. Yue, {\em Webster curvature and Hausdorff dimension of complex hyperbolic Kleinian groups}, In: ``Dynamical Systems --- Proceedings Of The International Conference In Honor Of Professor Liao Shantao,'' World Scientific, 1999. pp. 319--328.

\end{thebibliography}
\end{document}